\documentclass[3p]{scrartcl}%{article}%

\usepackage{bm, physics}
\usepackage{titling}
\usepackage[hmarginratio=1:1,top=32mm,left=15mm,columnsep=1cm]{geometry} % 
\usepackage{array}
\usepackage{color} 
\usepackage{tabularx}
\usepackage{graphicx}
\usepackage[svgnames,table,xcdraw]{xcolor}
\usepackage{amsmath}
\usepackage{amssymb}
\usepackage{amsthm}
\usepackage{amsfonts}	
\usepackage{comment}
\usepackage{multirow}
\usepackage{relsize} % e.g. used for \mathsmaller
\usepackage{fourier} % for the norm command \Vert (it has smaller space between the vertical bars 
\usepackage{epstopdf}
\usepackage[font=small,labelfont=bf]{caption}	% set font of caption text
\usepackage{subcaption}
\usepackage{multicol}
\usepackage{xcolor,cancel}
\usepackage{hyperref}
\usepackage{siunitx}

\hypersetup{
    colorlinks=true,                          
    linkcolor=Blue,
    citecolor=DarkRed,
    urlcolor=blue  }

\usepackage[backend=bibtex,giveninits=true,sorting=nyt,url=false,doi=true,eprint=false,isbn=false,
backref,backrefstyle=none,maxbibnames=99]{biblatex}
\DefineBibliographyStrings{english}{%
	backrefpage = {Cited on p\adddot},%
	backrefpages = {Cited on pp\adddot}%
}
\bibliography{references}

\renewbibmacro{in:}{%
	\ifboolexpr{%
		test {\ifentrytype{article}}%
	}{}{\printtext{\bibstring{in}\intitlepunct}}%
}
\newcommand\BibTeX{{\rmfamily B\kern-.05em \textsc{i\kern-.025em b}\kern-.08em
T\kern-.1667em\lower.7ex\hbox{E}\kern-.125emX}}

\newcommand*\samethanks[1][\value{footnote}]{\footnotemark[#1]}

\allowdisplaybreaks[1]

\DeclareMathOperator{\divg}{div}
\DeclareMathOperator{\diam}{diam}

\newcommand{\sgrad}[2]{  \langle \nabla #1  \rangle_{#2}}
\newcommand{\sdiv}[2]{  \langle \divg #1  \rangle_{#2}}
\newcommand{\adgrad}[2]{ \langle \nabla^* #1 \rangle_{#2}}
\newcommand{\addiv}[2]{  \langle \divg^* #1  \rangle_{#2}}

\newcommand{\vv}{\bm{v}}
\newcommand{\xx}{\bm{x}}
\newcommand{\yy}{\bm{y}}
\newcommand{\zz}{\bm{z}}

\newcommand{\dt}{\Delta t}

\renewcommand{\bm}{\boldsymbol}

\begin{document}
	
% TITLE
\title{ \Large
%	 \textbf{SHTC model for two-phase compressible flow with surface tension} 
	\textbf{A numerical method based on quasi-Lagrangian Voronoi cells for two-phase flows with large density contrast} 
}
% AUTHORS for article documentclass
\author{
Ondřej Kincl\thanks{Department of Civil, Environmental and Mechanical Engineering, 
University of Trento, Via Mesiano 77, Trento 38123, Italy
\href{mailto:ondrej.kincl@unitn.it}{ondrej.kincl@unitn.it}
},
\quad
Ilya Peshkov\samethanks[1],
\quad
Walter Boscheri\thanks{Laboratoire de Mathématiques UMR 5127 CNRS, Universit{\'e} Savoie Mont Blanc, 73376 Le Bourget du Lac, France }
}
\thanksmarkseries{arabic}
\maketitle %\maketitle must follow title, authors, abstract and \pacs
\date{\today}

% ABSTRACT
\begin{abstract}
\noindent
In this work, we use a moving Voronoi and sharp interface approach for simulating two-phase flows. At every time step, the mesh is generated anew from Voronoi seeds that behave as material points. The paper is a continuation of our previous works on moving Voronoi meshes where we have considered single phase incompressible and compressible flows. In the context of quasi-Lagrangian Voronoi simulations, problems with large density contrasts (such as water and air interface) are being treated here for the first time to the best of our knowledge. This is made possible through a remapping stage, which relies on a filtering of a color function. The resulting semi-implicit scheme is conservative and robust, allowing us to simulate both compressible and incompressible flows, including shock waves and surface tension. 
\end{abstract}

\section{Introduction}
In 1985, one of the the first Free-Lagrange method was proposed by Fritts, Crowley and Trease \cite{fritts1985free}. Their numerical scheme was based on a moving Voronoi tesselation, where generating seeds propagate with the flow and also serve as nodes for approximating derivatives in space. Voronoi cells are always convex and, regardless of the deformation, there is never any danger of producing inverted cells. As a result, the robustness can be compared to pure particle methods, like Smoothed Particle Hydrodynamics (SPH), but with certain benefits of cell-based approaches, like improved accuracy of gradient approximation and more compact stencils (i.e. each node communicates with a smaller number of neighbors at a given time). 

With the advent of more powerful computers, the potential of moving Voronoi grids was reinvestigated, for example, in the works of Gaburro et al \cite{gaburro2020high}, Springel \cite{springel2011hydrodynamic}, Després \cite{despres2024lagrangian}, and Loubère et al \cite{loubere2010reale}. In \cite{kincl2025semi} and \cite{kincl2025semi2}, we introduced a semi-implicit version of Free-Lagrange method, which we termed SILVA (Semi-Implicit Lagrangian Voronoi Approximation). The semi-implicit solver is advantageous for two reasons. Firstly, it avoids the stiffness problems associated with low-Mach flows. Secondly, by allowing larger time step, we also need fewer (relatively expensive) Voronoi mesh regenerations. Consequently, more time is spent on solving linear systems than mesh construction. SILVA was validated on both compressible and incompressible flows, featuring viscosity, gravity, shock discontinuities and multi-phase problems. But so far, the method was unstable in the case of high density contrast between the phases. This problem is resolved in the current paper. The improvement stems from a new, more robust, remapping technique, which uses a smoothed color function to stabilize the interface. This is an important achievement because, unlike in SPH, it is difficult to implement free surface problems in Free-Lagrange method. To the best of our knowledge, this has only been successfully done so far in the case of elasto-plastic solids by introducing virtual cells with zero material response \cite{howell2002free}.
 
This paper is organized as follows. In the methodology section, we explain some basic properties of the Voronoi grid and how derivatives are approximated. This formalism allows us to write a discrete compressible Navier-Stokes equation with surface tension. The time integration relies on a splitting technique, where reversible processes (pressure, surface tension and gravity) and irreversible forces (viscosity) are treated separately. Finally, to stabilize the mesh with respect to spurious clustering of generating seeds, a remapping stage based on Lloyd iterations is performed. In the results section, we validate the method by performing six numerical tests: the circular patch benchmark, dam breaking simulation, rotating square, rising bubble, shock-bubble interaction and interaction of a shock with a column of water.

The numerical code used for this paper is publicly available on GitHub \cite{Kincl_LagrangianVoronoi_jl}.

\section{Methodology}
    \subsection{Basic approximations}
\begin{figure}[hbt!]
        \centering
        \includegraphics[width=0.35\linewidth]{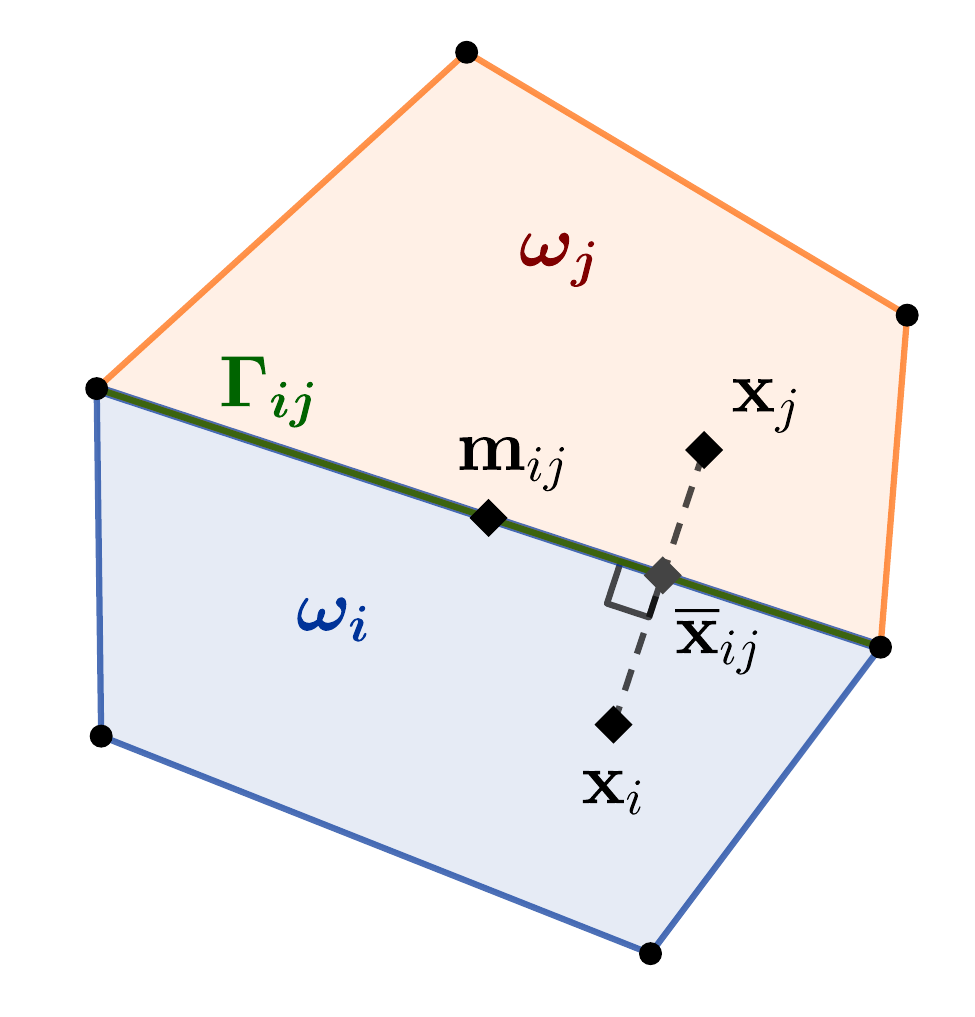}
        \caption{Two Voronoi cells in two-dimensional grid. Note that the midpoint between the generating seeds $\overline{\xx}_{ij} = \frac{1}{2}(\xx_i + \xx_j)$ is generally distinct from the edge midpoint $\bm{m}_{ij}$.}
        \label{fig:cells}
    \end{figure}

    Let the computational domain $\Omega \subset \mathbb{R}^d$ be a convex polytope, where $d$ is either 2 or 3. Let $\xx_i \in \Omega$ for $i = 1,\dots,N$ be a collection of \textit{distinct} points. These points serve as computational nodes and are also the generating seeds of the Voronoi mesh. The $i$-th Voronoi cell is defined as follows:
    \begin{equation}
        \omega_i = \bigcap_{j \neq i} \{\xx \in \Omega : |\xx - \xx_j| < |\xx - \xx_i|\}.
    \end{equation}
    The Voronoi cell $\omega_i$ is also a convex polytope. Two Voronoi cells $\omega_i$ and $\omega_j$ are neighbors if 
    \begin{equation}
        \Gamma_{ij} = \partial \omega_i \cap \partial \omega_j
    \end{equation}
    is non-empty. Let $N(i)$ be the set of all indices $j$ such that $\omega_i$ and $\omega_j$ are neighbors. An example of neighboring Voronoi cells is depicted in Figure \ref{fig:cells}. 
    For $j \in N(i)$, the surface normal of facet $\Gamma_{ij}$ pointing from $\omega_i$ to $\omega_j$ can be expressed by taking the difference of the seeds:
    \begin{equation}
        \bm{n}_{ij} = \frac{\xx_j - \xx_i}{|\xx_j - \xx_i|}.
    \end{equation}
    This is true by virtue of the orthogonality between the equidistant plane (line, for $d=2$)
    \begin{equation}
        P_{ij} = \bigg\{ \xx \in \mathbb{R}^d:|\xx-\xx_i| = |\xx-\xx_j| \bigg\} \supseteq \Gamma_{ij}
    \end{equation}
    and the vector of separation $\xx_i - \xx_j$. The tesselation allows to approximate integrals by
    \begin{equation}
        \int_\Omega f(\xx) \, \dd \xx \approx \sum_i |\omega_i| f_i,
    \end{equation}
    where $f_i = f(\xx_i)$ and $|\omega_i|$ denotes the $d$-dimensional measure of cell $\omega_i$. Using rough estimates, it is easy to infer that
    \begin{equation}
        \left| \int_\Omega f(\xx) \, \dd \xx - \sum_i |\omega_i| f_i \right| \leq |\Omega| \| \nabla f\|_\infty \max_i h_i,
    \end{equation}
    where $h_i = \diam \omega_i$ is the diameter of the cell $\omega_i$, provided that $f$ is continuously differentiable on ${\Omega}$. Furthermore, the gradient of a scalar field $f$ can be conveniently approximated by
    \begin{equation}
        \langle \nabla f\rangle_i =- \frac{1}{|\omega_i|} \sum_{j \in N(i)} \frac{|\Gamma_{ij}|}{r_{ij}}(f_i - f_j)(\bm{m}_{ij} - \xx_i).
        \label{eq:sgrad}
    \end{equation}
    Here, $|\Gamma_{ij}|$ is the $(d-1)$-dimensional measure of the facet, $r_{ij} = |\xx_i - \xx_j|$ and 
    \begin{equation}
        m_{ij} = \int_{\Gamma_{ij}} \xx \, \dd S(\xx)
    \end{equation}
    is its centroid. It can be proven that such approximation is first order consistent. That means
    \begin{equation}
        |\langle \nabla f\rangle_i - \nabla f(\xx_i)| = O(h_i),
    \end{equation}
    given that $f$ is twice continuously differentiable on ${\omega}_i$ and provided that either $\omega_i$ is an interior cell ($\partial \omega_i \cap \partial \Omega = \emptyset$) or $f$ satisfies a homogeneous Neumann condition on $\partial \omega_i \cap \partial \Omega$. The proof can be found in \cite{springel2010pur}. Additionally, there exists a dual approximation
    \begin{equation}
        \begin{split}
            \langle \nabla^* f \rangle_i =& \frac{1}{|\omega_i|} \sum_{j\in N(i)} \frac{|\Gamma_{ij}|}{r_{ij}}\bigg[f_i(\bm{m}_{ij} -\xx_i) - f_j(\bm{m}_{ij} - \xx_j) \bigg]\\
            =&\frac{1}{|\omega_i|} \sum_{j\in N(i)} \frac{|\Gamma_{ij}|}{r_{ij}} \bigg( (f_{i}-f_j)\left(\bm{m}_{ij} - \frac{1}{2}(\xx_i+\xx_j)\right) -  \frac{1}{2}(f_i+f_j) (\xx_{i}-\xx_j) \bigg),
        \end{split}
        \label{eq:wgrad}
    \end{equation}
    which can be motivated by a discrete analogy of integration by parts:
    \begin{equation}
        \sum_i |\omega_i| \langle \nabla f_i\rangle  g_i = - \sum_i |\omega_i| f_i \langle \nabla^* g\rangle_i.
        \label{eq:byparts}
    \end{equation}
    Unlike \eqref{eq:sgrad}, such approximation is consistent only in a weak sense and is usually a bit noisy \cite{kincl2025semi}. Nonetheless, this operator is needed to obtain a conservative scheme. Any first order partial derivative can be approximated using \eqref{eq:wgrad} or \eqref{eq:sgrad} by extracting a specific component of the gradient. In other words, the partial derivative of $f$ with respect to $x,y$ and $z$ coordinates is obtained as follows:
    \begin{align}
        \left \langle \pdv{f}{x} \right \rangle_i &= \hat{\xx}_i \cdot \langle \nabla f_i\rangle\\
        \left \langle \pdv{f}{y} \right \rangle_i &= \hat{\yy}_i \cdot \langle \nabla f_i\rangle\\
        \left \langle \pdv{f}{z} \right \rangle_i &= \hat{\zz}_i \cdot \langle \nabla f_i\rangle,
    \end{align}
    where $\hat{\xx}, \hat{\yy}$ and $\hat{\zz}$ are the basis of a Cartesian coordinate system. This allows to approximate the divergence of a vector field by
    \begin{equation}
        \langle \divg \bm{v}\rangle_i = \left \langle \pdv{v^x}{x} \right \rangle + \left \langle \pdv{v^y}{y} \right \rangle + \left \langle \pdv{v^z}{z} \right \rangle =-\frac{1}{|\omega_i|} \sum_{j \in N(i)} \frac{|\Gamma_{ij}|}{r_{ij}}(\vv_i - \vv_j) \cdot (\bm{m}_{ij} - \xx_i).
        \label{eq:sdiv}
    \end{equation}
    This algebraic compatibility is used everywhere throughout this paper and applies to all first order differential operators such as curl or divergence of a tensor field. Approximation based on the dual gradient \eqref{eq:wgrad} will be always denoted by a star. For instance,
    \begin{equation}
        \begin{split}
             \langle \divg^* \bm{v} \rangle_i &= \frac{1}{|\omega_i|} \sum_{j\in N(i)} \frac{|\Gamma_{ij}|}{r_{ij}}\bigg[\vv_i \cdot (\bm{m}_{ij} -\xx_i) - \vv_j\cdot(\bm{m}_{ij} - \xx_j) \bigg].\\
        \end{split}
        \label{eq:wdiv}
    \end{equation}
    An interesting identity arises, when we evaluate the rate of change of cell volume with respect to moving generating seeds. It reads formally
    \begin{equation}
        \dv{|\omega_i|}{t} = |\omega_i| \left \langle \divg^* \bm{u} \right \rangle_i,
        \label{eq:dvolume}
    \end{equation}
    where $\bm{u}_i = \dv{\xx_i}{t}$ denotes the seed velocity for each $i$. This is a corollary of Reynolds transport theorem. See \cite{despres2024lagrangian} for a detailed proof.
    
    \subsection{Two-phase fluid flows}
    Firstly, we shall study a two-phase flow in homogeneous gravitational field with surface tension, where each component can be described as a compressible Navier-Stokes fluid. Using the material time derivative
	\begin{equation}
		\dv{f}{t} = \pdv{f}{t} + \vv \cdot \nabla f, \quad \forall f,
	\end{equation}
    the system of partial differential equations in consideration can be written in the following Lagrangian form:
    \begin{align}
		\dv{\rho}{t} =&  -\rho \divg{\vv}, \label{eq:bomass}\\
		\rho \dv{\vv}{t}  =&  \divg{\mathbb{T}} + \bm{g}, \label{eq:bomom}\\
		\rho \dv{e}{t}  =& (\divg{\mathbb{T}})\cdot \bm{v} + \mathbb{T}:\nabla \vv , \label{eq:boe}\\
        \dv{C}{t} =& 0, \label{eq:C}
	\end{align}
    where $\rho$ is the fluid density, $\vv$ is the fluid velocity, $e$ is the total specific energy, $\mathbb{T}$ is the Cauchy stress tensor, $\bm{g}$ is gravitational acceleration and, lastly, $C$ is the color function. We restrict ourselves to a sharp interface model between two phases, so that the variable $C$ can only take two different values: 
    \begin{equation}
        C(t,\xx) = \begin{cases}
         0, & \xx\in\text{primary phase},\\
         1, & \xx\in\text{secondary phase}.
        \end{cases}
    \end{equation}
    The stress tensor is
    \begin{equation}
        \mathbb{T} = -p\mathbb{I} + 2\mu \mathbb{D} + \mathbb{S},
    \end{equation}
    where $p = p(\rho,\epsilon,C)$ is the pressure specified by an equation of
    state (different for each phase in general), $\mathbb{D} =
    \frac{1}{2}(\nabla \vv + \nabla \vv^T)$ is the symmetric part of the
    strain-rate tensor and the parameter $\mu = \mu(C)$ is the dynamic viscosity. To
    determine $\epsilon$, the specific internal energy, we need to subtract
    kinetic and tensile components:
    \begin{equation}
        \epsilon = e - \frac{\vv^2}{2} + \bm{g}\cdot \xx - \sigma |\nabla C^H|,
    \end{equation}
    where the gravitational potential is also considered. The last component $\mathbb{S}$ of the stress tensor stands for the stress by surface tension. Following \cite{violeau}, we assume it has the form: 
    \begin{equation}
        \mathbb{S} = \sigma |\nabla C^H|\bigg(\mathbb{I} - \bm{n}^I \otimes \bm{n}^I\bigg),
    \end{equation}
    where $\sigma$ is the surface tension coefficient,
    \begin{equation}
        \bm{n}^I = \frac{\nabla C^H}{|\nabla C^H|}
    \end{equation}
    is the interface normal and 
    \begin{equation}
        C^H(\xx) = (w^H*C)(\xx) = \int C(\yy) w^H(\xx - \yy)\dd \yy
    \end{equation}
    is a color function convoluted with a Wendland's quintic kernel 
    \begin{equation}
        w_H(\xx) =
					\begin{cases}
						\frac{\alpha_d}{H^2} \left( 1 - \frac{r}{H} \right)^4 \left(1 +  \frac{4r}{H}\right), & \mathrm{for} \quad r < H\\
						0, & \mathrm{for} \quad r \geq H
					\end{cases}
    \end{equation}
    where the smoothing radius $H$ determines the width of the support for tensile forces and
    \begin{equation}
        \alpha_d = \begin{cases}
            \frac{3}{2}, & d = 1\\
            \frac{7}{\pi}, & d=2\\
            \frac{21}{2\pi}, & d=3
        \end{cases}
    \end{equation}
    is a dimension-dependent normalization constant.
    Using the generating seeds as the simulation nodes, we can directly translate the governing equations \eqref{eq:bomass}-\eqref{eq:C} into a semi-discrete approximation (a system of ordinary differential equations where time remains continuous):
    \begin{align}
        \dv{\xx_i}{t} =\,& \vv_i\\
        \dv{\rho_i}{t} =\,&  -\rho_i \addiv{\vv}{i}, \\
		\rho_i \dv{\vv_i}{t}  =\,&  \sdiv{\mathbb{T}}{i} + \bm{g}, \\
		\rho_i \dv{e_i}{t}  =\,& \sdiv{\mathbb{T}}{i}\cdot \bm{v}_i + \mathbb{T}_i:\adgrad{\vv}{i}, \\
        \dv{C_i}{t} =\,& 0,
        \label{eq:twophase}
    \end{align}
    where
    \begin{gather}
        \mathbb{T}_i = -p_i\mathbb{I} + 2\mu_i \mathbb{D}_i + \mathbb{S}_i,\\[2mm]
        \mathbb{D}_i = \frac{1}{2}\bigg(\adgrad{\vv}{i} + \adgrad{\vv}{i}^T\bigg),  \\[2mm]
        \mathbb{S}_i = \sigma |(\nabla C^H)_i|\left( 
        \mathbb{I} - \frac{(\nabla C^H)_i  \otimes (\nabla C^H)_i}{|(\nabla C^H)_i|^2} \right).
        \label{eq:def_of_D_i}
    \end{gather}
     The gradient of a color function can be approximated like in SPH:
    \begin{equation}
        \nabla C^H(\xx_i) = \int \nabla w^H(\xx_i-\yy)\, C(\yy)\dd \yy \approx \sum_i |\omega_i| C_j \nabla w^H(\xx_i-\xx_j) =: (\nabla C^H)_i,
    \end{equation}
    where $H$ is linked to the spatial step by $H = 3h$, implying that the smoothing is made over the length of three cells.
    Our combination of starred and unstarred operators is not arbitrary and can be deduced from Euler-Lagrange equations, which we will not detail here (we point the interested reader to \cite{serrano2005voronoi}). Let us only mention that the discrete continuity equation together with relation \eqref{eq:dvolume} guarantees conservation of the cell mass 
    \begin{equation}
        M_i = \rho_i |\omega_i|. \label{eq:mass_identity}
    \end{equation}
    Verily, we have
    \begin{equation}
        \dv{M_i}{t} = \dv{\rho_i}{t} |\omega_i| + \rho_i \dv{|\omega_i|}{t} = -\rho_i |\omega_i| \addiv{\vv}{i} + \rho_i |\omega_i| \addiv{\vv}{i} = 0.
    \end{equation}
    The property \eqref{eq:mass_identity} is useful, because it allows to compute density directly from the initial mass. The total energy is conserved by virtue of the discrete integration by parts \eqref{eq:byparts}. Assuming periodic boundary condition or the absence of stress near boundaries, the linear momentum is also provably conserved, since
    \begin{equation}
        \sum_i |\omega_i| \rho_i \dv{\vv_i}{t} = \sum_i |\omega_i|\sdiv{\mathbb{T}}{i} = - \sum_i |\omega_i|{\mathbb{T}}_{i}\adgrad{1}{i},
    \end{equation}
    where the following identity holds for every interior cell:
    \begin{equation}
        \adgrad{1}{i} = \frac{1}{|\omega_i|} \sum_{j\in N(i)} \frac{|\Gamma_{ij}|}{|\xx_i - \xx_j|}\bigg[ (\bm{m}_{ij} -\xx_i) -(\bm{m}_{ij} - \xx_j) \bigg] = \frac{1}{|\omega_i|} \sum_{j\in N(i)}|\Gamma_{ij}|\bm{n}_{ij} = \frac{1}{|\omega_i|}\int_{\partial \omega_i} \bm{n} \, \dd S = 0.
    \end{equation}

    \subsection{Time-marching scheme}
    Our time integration relies on operator splitting and is similar to the work forwarded in \cite{kincl2025semi2}, except that we now include surface tension.

    Let us denote the numerical solution at time $t^{(n)}$ as 
    \begin{equation}
        Q^{(n)} = \begin{pmatrix}
            \xx_i^{(n)}\\
            \rho_i^{(n)}\\
            \vv_i^{(n)}\\
            e_i^{(n)}\\
            C_i^{(n)}
        \end{pmatrix}.
    \end{equation}
In every time step $\dt=t^{n+1}-t^n$, the variables are updated as follows: 
    \begin{equation}
        Q^{(n+1)} = \varphi^\mathrm{rev}_{\dt}( \varphi^\mathrm{irr}_{\dt}(Q^{(n)}))
    \end{equation}
    where $\varphi^\mathrm{rev}$ is an approximation of a flow operator corresponding to the \textit{reversible} part of the dynamics:
    \begin{align}
        \dv{\xx_i}{t} =\,& \vv_i \label{eq:path_rev}\\
        \dv{\rho_i}{t} =\,&  -\rho_i \addiv{\vv}{i}, \label{eq:bomass_rev}\\
		\rho_i \dv{\vv_i}{t}  =\,&  \sdiv{\mathbb{T}^\mathrm{rev}}{i} + \bm{g}, \label{eq:bomom_rev}\\
		\rho_i \dv{e_i}{t}  =\,& \sdiv{\mathbb{T}^\mathrm{rev}}{i}\cdot \bm{v}_i + \mathbb{T}^\mathrm{rev}_i:\adgrad{\vv}{i}, \label{eq:boe_rev} \\
        \dv{C_i}{t} =\,& 0, \label{eq:color_rev}
    \end{align}
    whereas $\varphi^\mathrm{irr}$ is concerned with the \textit{irreversible} (viscous) dynamics:
        \begin{align}
            \dv{\xx_i}{t} =\,& \dv{\rho_i}{t} = \dv{C_i}{t} = 0, \label{eq:path_irr}\\
    		\rho_i \dv{\vv_i}{t}  =\,&  \sdiv{\mathbb{T}^\mathrm{irr}}{i}, \label{eq:bomass_irr}\\
    		\rho_i \dv{e_i}{t}  =\,& \sdiv{\mathbb{T}^\mathrm{irr}}{i}\cdot \bm{v}_i + \mathbb{T}^\mathrm{irr}_i:\adgrad{\vv}{i}, \label{eq:boe_irr} \\
    \end{align}
    where
    \begin{equation}
        \mathbb{T}^\mathrm{rev} =-p_i\mathbb{I} + \mathbb{S}_i,
        \qquad
        \mathbb{T}^\mathrm{irr} =2\mu_i \mathbb{D}_i.
    \end{equation}
    The flow operators $\varphi^\mathrm{rev}_{\dt}$ and $\varphi^\mathrm{irr}_{\dt}$ are constructed in the following subsections.
    
        \subsection{Irreversible step}
    This section is devoted to the construction of the irreversible fractional step 
    \begin{equation}
        \varphi^\mathrm{rev}_{\dt}:Q^{(n)} \to Q^{(n+1)}.
    \end{equation}
    There are two options: an implicit and an explicit one. The explicit
    approach is computationally cheap, easy to implement and is sufficiently
    robust. However, it can severely under-perform in case of flows with very
    low Reynolds number. To simplify notation, we describe the
    explicit step for $n=0$ with no loss of generality:
    \begin{align}
		\xx_i^{(1)} &= \xx_i^{(0)}, \\
		\rho_i^{(1)} &= \frac{M_i}{|\omega_i^{(1)}|}, \\
		\frac{\vv_i^{(1)} - \vv_i^{(0)}}{\dt}  &=  \frac{2}{\rho_i^{(0)}}  \sdiv{(\mu \mathbb{D}^{(0)})}{i}, \\
		\frac{e_i^{(1)} - e_i^{(0)}}{\dt}  &=  \frac{2}{\rho_i^{(0)}}\left(\sdiv{(\mu \mathbb{D}^{(0)})}{i}\cdot \vv_i^{(0)} + (\mu_i \mathbb{D}^{(0)}_i):\adgrad{\vv^{(0)}}{i}\right),  \\
        C_i^{(1)} &= C_i^{(0)}.
	\end{align}
    Regarding the implicit option, with $n=1$ referring to the next time step, it reads as follows:
        \begin{align}
		\xx_i^{(1)} &= \xx_i^{(0)}, \\
		\rho_i^{(1)} &= \frac{M_i}{|\omega_i^{(1)}|}, \\
		\frac{\vv_i^{(1)} - \vv_i^{(0)}}{\dt}  &=  \frac{2}{\rho_i^{(1)}}  \sdiv{(\mu \mathbb{D}^{(1)})}{i}, \\
		\frac{e_i^{(1)} - e_i^{(0)}}{\dt}  &=  \frac{2}{\rho_i^{(1)}}\left(\sdiv{(\mu \mathbb{D}^{(1)})}{i}\cdot \vv_i^{(1)} + (\mu_i \mathbb{D}^{(1)}_i):\adgrad{\vv^{(1)}}{i}\right),  \\
        C_i^{(1)} &= C_i^{(0)}.
	\end{align}
    It is then necessary to solve the linear system:
    \begin{equation}
        \frac{M_i\vv_i^{(1)}}{\dt}-  2|\omega_i|\sdiv{(\mu \mathbb{D}^{(1)})}{i} = \frac{M_i \vv_i^{(0)}}{\dt}
        \label{eq:bomom_irr2} 
    \end{equation}
    for the unknown velocity field $\vv^{(1)}$. The system \eqref{eq:bomom_irr2} is positive definite, which we can verify by using \eqref{eq:byparts} and \eqref{eq:def_of_D_i} to calculate
    \begin{equation}
        \begin{split}
            \sum_i\left(\frac{M_i\vv_i}{\dt}-  2|\omega_i|\sdiv{(\mu \mathbb{D})}{i}\right) \cdot \vv_i &= \sum_i\frac{M_i|\vv_i|^2}{\dt} + 2\sum_i |\omega_i| \mu_i \mathbb{D}:\nabla \vv_i\\
            &= \sum_i\frac{M_i|\vv_i|^2}{\dt} + 2\sum_i |\omega_i| \mu_i| \mathbb{D}|^2 > 0,
        \end{split}
    \end{equation}
    for arbitrary velocity field $\vv$. Unfortunately, the associated matrix is not very sparse, because an expensive double summation is involved in $\sdiv{(\mu \mathbb{D})}{}$. This limitation can be avoided with a matrix-free conjugate method, where the velocity gradient (the $\mathbb{D}$-field) is saved into memory temporarily during each application of the linear operator. Once $\vv^{(1)}$ is known, updating the energy becomes an explicit process. We employ this implicit scheme in the rising bubble test, where the Reynolds number is relatively small.
    
    \subsection{Reversible step}
    This subsection describes the map
    \begin{equation}
        \varphi^\mathrm{rev}_{\dt}:Q^{(n)} \mapsto Q^{(n+1)},
    \end{equation}
    approximating the advance of time by $\dt$ if viscosity is absent. To simplify the notation, let us assume $n=0$. Because of the stiff relation between density and pressure in low-Mach regime, we require the scheme to be semi-implicit. The updated values are defined as follows:
    \begin{align}
		\frac{\xx_i^{(1)} - \xx_i^{(0)}}{\dt} &= \vv_i^{(0)}, \label{eq:path_implicit}\\
		\rho_i^{(1)} &= \frac{M_i}{|\omega_i^{(1)}|}, \label{eq:bomass_implicit}\\
		\frac{\vv_i^{(1)} - \vv_i^{(0)}}{\dt}  &=  \frac{1}{\rho_i^{(1)}} \left( -\sgrad{\tilde{p}}{i} + \sdiv{\mathbb{S}^{(1)}}{i}\right)+ \bm{g}, \label{eq:bomom_implicit}\\
		\frac{e_i^{(1)} - e_i^{(0)}}{\dt}  &= \frac{1}{\rho_i^{(1)}}\left( -\sgrad{\tilde{p}}{i} \cdot \vv_i^{(1)} - \tilde{p}_i \addiv{\vv^{(1)}}{i} + \sdiv{\mathbb{S}^{(1)}}{i}\cdot\vv_i^{(1)} \right),  \label{eq:boe_implicit}\\
        C_i^{(1)} &= C_i^{(0)}.
	\end{align}
    Here, all differential operators are computed with respect to the updated positions $\xx^{(1)}$, and $\tilde{p}$ is a certain prediction of $p^{(1)}$ that can be computed easily (without having to solve a non-linear problem). We start by writing
    \begin{equation}
        \dv{p}{t} = \left(\pdv{p}{\rho}\right)_{\eta} \dv{\rho}{t} = c^2 \dv{\rho}{t}.
    \end{equation}
    This holds because the specific entropy $\eta$ is preserved by the reversible part of the dynamics. Without surface tension, i.e. $\sigma = 0$, this property also holds for the semi-discrete system \eqref{eq:bomass_rev}-\eqref{eq:boe_rev}. The proof is exactly the same as in \cite{kincl2025semi2}. It could be possible to derive an entropic robust and variationally consistent implementation of surface tension by writing the Euler-Lagrangian equation using an approximation of tensile energy. This however leads to great complications. Using \eqref{eq:bomass} leads to:
    \begin{equation}
        \frac{\tilde{p} - p_i^{(0)}}{\dt \, \left(c_i^{(0)}\right)^2} = -\rho_i^{(1)}  \addiv{\vv^{(1)}}{i}= -\rho_i^{(1)}  \addiv{\tilde{\vv}}{i}+ \dt \rho_i^{(1)} \left \langle {\divg^*\frac{\sgrad{\tilde{p}}{}}{\rho^{(1)}}}\right \rangle_{i}\label{eq:tilde_p_implicit},
    \end{equation}
    where
    \begin{equation}
        \tilde{\vv} = \vv^{(0)}_i + \frac{\dt}{\rho_i^{(1)}}\sdiv{\mathbb{S}^{(1)}}{i} + \bm{g}.
    \end{equation}
    Note that surface tension does not depend on velocity or energy and so $\tilde{\vv}$ can be computed explicitly. In order to narrow the stencil and avoid problems with odd-even decoupling \cite{date2003fluid}, the composed operator in \eqref{eq:tilde_p_implicit} is further approximated using
    \begin{equation}
		\sgrad{\tilde{p}}{i}\cdot (\xx_i^{(1)}-\xx_j^{(1)}) \approx \tilde{p}_i - \tilde{p}_j \approx \sgrad{\tilde{p}}{j} \cdot (\xx_i^{(1)} - \xx_j^{(1)}),
	\end{equation}
    which, upon expanding \eqref{eq:wdiv}, leads to
	\begin{equation}
		\begin{split}
			-|\omega_i|^{(1)}\bigg \langle {\divg}^* { \frac{\langle \nabla \tilde{p}\rangle}{\rho^{(1)}} }\bigg \rangle_{i} &= \sum_{j\in N(i)}\frac{|\Gamma_{ij}^{(1)}|}{r_{ij}^{(1)}} \left(   \left( \frac{\sgrad{\tilde{p}}{i}}{2\rho_i^{(1)}} + \frac{\sgrad{\tilde{p}}{j}}{2\rho_j^{(1)}} \right)\cdot (\xx_{i}^{(1)}-\xx_j^{(1)})- \left( \frac{\sgrad{\tilde{p}}{i}}{\rho_i^{(1)}} - \frac{\sgrad{\tilde{p}}{j}}{\rho_j^{(1)}} \right)\cdot \left(\bm{m}_{ij}^{(1)} - \frac{1}{2} (\xx_i + \xx_j)^{(1)}\right) \right) \\
			&\approx \sum_{j\in N(i)} \frac{|\Gamma_{ij}^{(1)}|}{r_{ij}^{(1)}} \left(   \left( \frac{1}{2\rho_i^{(1)}} + \frac{1}{2\rho_j^{(1)}} \right)(\tilde{p}_{i}-\tilde{p}_j)- \left( \frac{\sgrad{p^\dagger}{i}}{\rho_i^{(1)}} - \frac{\sgrad{p^\dagger}{j}}{\rho_j^{(1)}} \right)\cdot \left(\bm{m}_{ij}^{(1)} - \frac{1}{2}(\xx_i + \xx_j)^{(1)}\right) \right),\end{split}
		\label{eq:sparsification_trick}
	\end{equation}
    where $p^\dagger$ in the second term is the last known pressure value, starting with $p^\dagger = p^{(0)}$ and then proceeding via fixed point iterations. In each iteration, the resulting system features a sparse positive semi-definite matrix and can be solved by conjugate gradient method. For details, we refer to our previous paper \cite{kincl2025semi2}. This technique copes with low-Mach flows and even incompressible flows ($c = \infty$). Once the pressure $\tilde{p}$ is known, it can be substituted into \eqref{eq:bomom_implicit} in order to obtain the updated velocity, which is finally used to increment the energy in \eqref{eq:boe_implicit}.
    
    \subsection{Remapping}
    Attempts to make a purely Lagrangian method based on Voronoi tessellations usually fail because of spurious particle clustering. A Voronoi mesh is ill defined when the generating seeds overlap. However, in finite arithmetic, even a small inter-seed distance is problematic. A convenient dimensionless measure for the cell quality is the minimum ratio between seed distances. That is,
    \begin{equation}
        q_\mathrm{cell,i} = \frac{\min_{j\in N(i)}r_{ij}}{\max_{k\in N(i)}r_{ik}}.
    \end{equation}
    The maximum value is $q_\mathrm{cell,i} = 1$ and can be found in regular Cartesian or hexagonal lattices. The overall mesh quality can be defined as
    \begin{equation}
        q_\mathrm{mesh} = \min_i q_\mathrm{cell,i}.
    \end{equation}
    In our scheme, the mesh is repaired when $q_\mathrm{mesh}$ drops below
    $q_\mathrm{treshold} = 0.3$ by a single Lloyd iteration. This means that the
    generating seeds are redefined to be the centroids of the respective cell.
    To prevent artificial mixing between different phases, the centroid is
    weighted by a smoothed color function. Mathematically:
    \begin{equation}
        \xx_i^\mathrm{remap} = \begin{cases}
            \frac{\int_{\omega_i} C^H(\xx) \xx  \, \dd \xx}{\int_{\omega_i} C^H(\xx) \, \dd \xx} & \mathrm{for} \;C_i = 1,\\
            \\
            \frac{\int_{\omega_i} (1-C^H(\xx)) \xx \, \dd \xx}{\int_{\omega_i} (1-C^H(\xx)) \, \dd \xx} & \mathrm{for} \; C_i = 0.
        \end{cases}
        \label{eq:centroid}
    \end{equation}
    Integrals in \eqref{eq:centroid} are computed numerically by splitting the Voronoi cell into simplexes and using a second order quadrature rule on every simplex. To restore consistency, the variables need to be remapped by explicitly solving a linear advection problem. To remain conservative, the update must be formulated in terms of the cell mass 
    \begin{equation}
        M = \rho_i|\omega_i|,
    \end{equation}
    momentum
    \begin{equation}
        \bm{U}_i = \rho_i |\omega_i| \vv_i
    \end{equation}
    and total energy
    \begin{equation}
        E_i = \rho_i |\omega_i| e_i.
    \end{equation}
    The update from old variables to the remapped variables is:
    \begin{align}
		M_i^\mathrm{remap} - M_i &= \mathcal{F}_i(\rho), \label{eq:remap_M} \\
		\bm{U}_i^\mathrm{remap} - \bm{U}_i &= \mathcal{F}_i{(\rho \vv )}, \label{eq:remap_U} \\
		E_i^\mathrm{remap} - E_i &= \mathcal{F}_i{(\rho e)},\label{eq:remap_E}
	\end{align}
    where the Rusanov flux $\mathcal{F}_i$ is defined as
    \begin{equation}
            \mathcal{F}_i(\phi) =  \sum_{\begin{matrix} j \in N(i) \\ C_i = C_j\end{matrix}} \frac{|\Gamma_{ij}|}{r_{ij}}\bigg[ \varphi _i\delta\xx_i \cdot (\bm{m}_{ij} -\xx_i) - \varphi_j \delta \xx_j\cdot(\bm{m}_{ij} - \xx_j) -\frac{r_{ij}}{2}\max \{ |\delta \bm{x}_i|, |\delta \bm{x}_j| \}\;(\phi_i - \phi_j )\bigg], \quad \forall \phi,
    \end{equation}
    and $\delta \xx_i = \xx_i^\mathrm{remap} - \xx_i$. Note that the summation must be restricted to cells of the same phase. Otherwise, an unphysical exchange of mass, momentum and energy occurs between different phases and the color function $C$ would also need to be remapped, making a previously sharp interface diffuse. Unfortunately, it is currently not clear if such scheme remains consistent near the interface. Nonetheless, our numerical experiments indicate linear order of convergence.

\begin{figure}[htb!]
        \centering
        \includegraphics[width=0.5\linewidth]{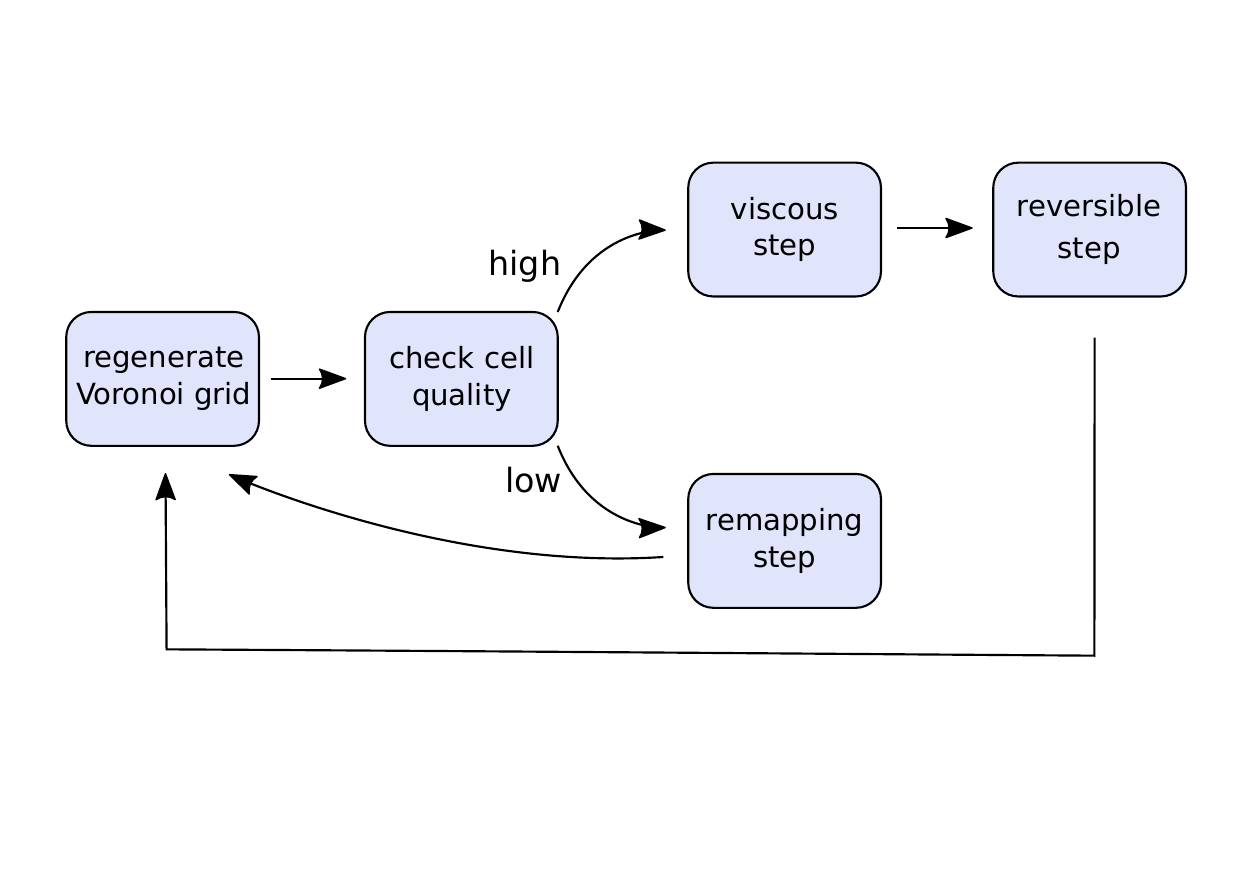}
        \caption{Flowchart explaining the time-marching scheme of SILVA.}
        \label{fig:flowchart}
    \end{figure}

    The flowchart in Figure \ref{fig:flowchart} shows the order of steps in our scheme.

\section{Numerical results}
    \subsection{Circular patch}
    There exists a simple free surface problem, called circular patch, which
    admits a semi-analytical solution \cite{monaghan2000sph}. Unlike in SPH,
    free surface is challenging to implement in a Voronoi method, since the shape of
    the fluid needs to be known \textit{a priori} to generate Voronoi cells and it is
    difficult to infer from generating seeds alone. Therefore, our scheme cannot
    handle a free surface problem directly. However, we can mimic free surface
    behavior via two-phase flow with a large density contrast. 
    
    In the circular patch problem, the computational domain is a rectangle 
    $\Omega = \left(-3/2,3/2\right) \times (-3,3) $
    and contains two incompressible fluids. In the initial frame, water occupies a unit circle and has a prescribed velocity field
    \begin{equation}
    \vv = \frac{1}{2}\begin{pmatrix}
        -x\\
        y
    \end{pmatrix}.
    \end{equation}
    Air, initially at rest, occupies remainder of the domain. We prescribe free-slip boundaries everywhere. The water density is $\rho_\mathrm{water} = 1000 \unit{Kg/m^3}$ and the air density is $\rho_\mathrm{air} = 1.25 \unit{Kg/m^3}$, so that the density contrast is 
    $ \frac{\rho_\mathrm{water}}{\rho_\mathrm{air}} = 800$.
    \begin{figure}[htb!]
        \centering
        \begin{subfigure}{0.5\linewidth}
            \centering
            \includegraphics[width=\linewidth]{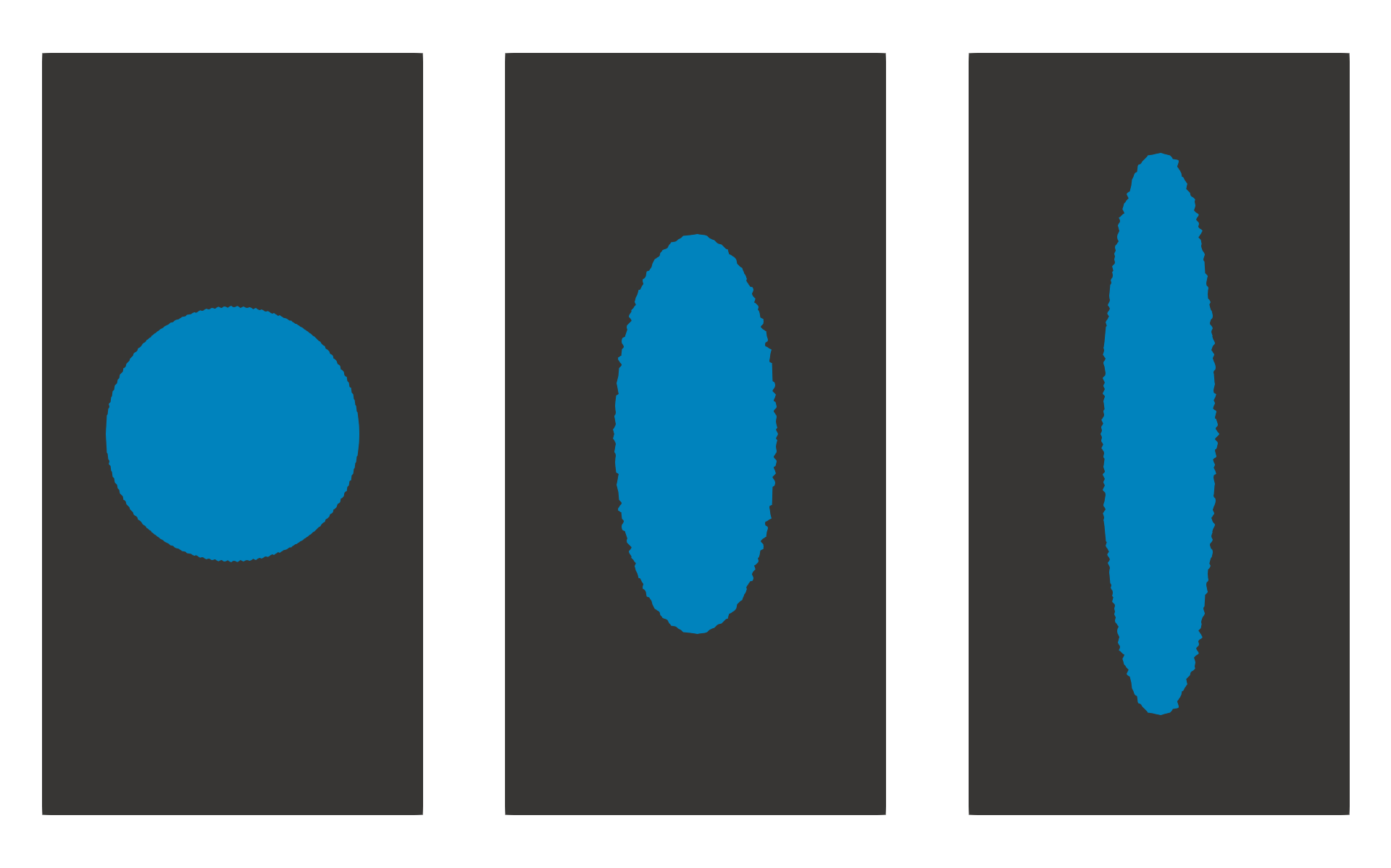}
            \caption{The color function at $t = 0$ (left), $t = 1$ (middle) and $t=2$ (right).}
        \end{subfigure}%
        \begin{subfigure}{0.5\linewidth}
            \centering
            \includegraphics[width=\linewidth]{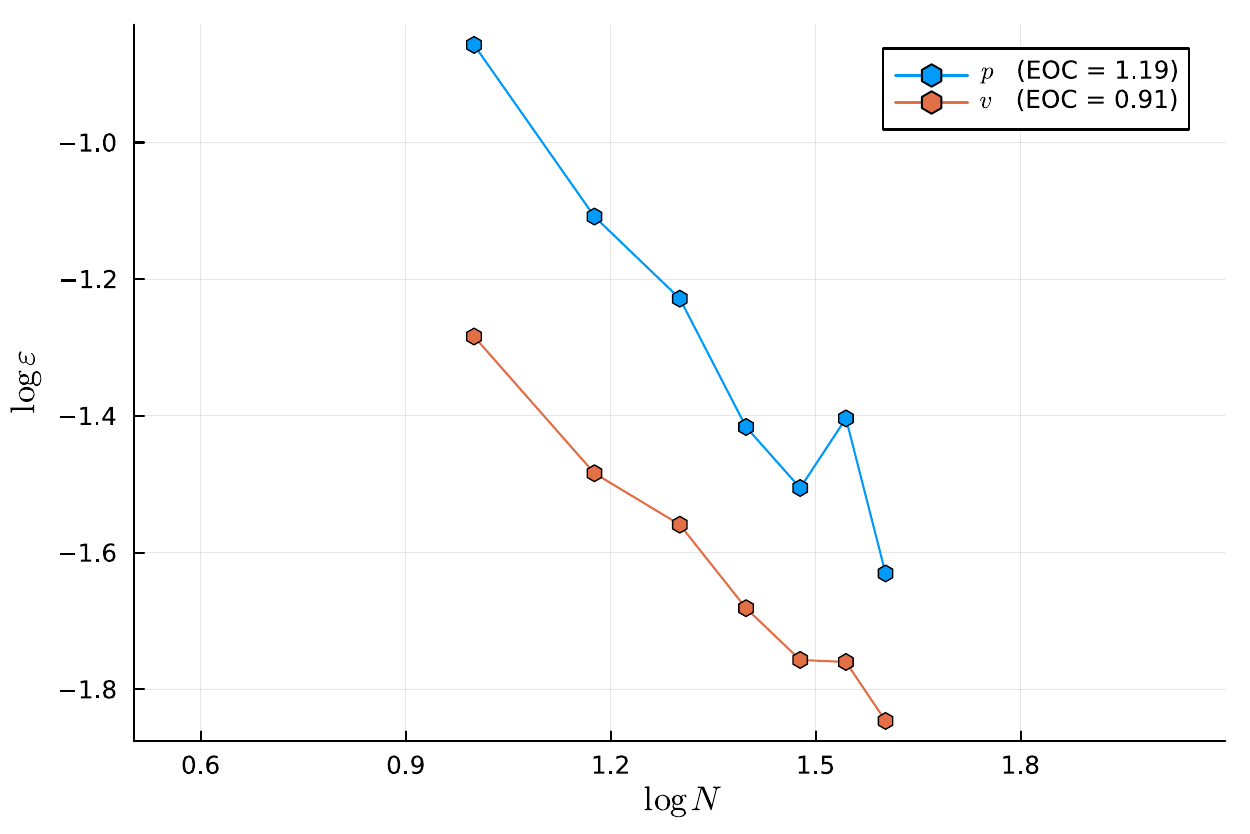}
            \caption{Convergence study at $t=3 \unit{s}$. Here, $N$ is the number of sample points per radius of the circle and $\epsilon$ is error in the $L^2$ norm.}
        \end{subfigure}
        \caption{Circular path test results.}
        \label{fig:cpatch}
    \end{figure}
    The effects of viscosity and surface tension are neglected. As precluded earlier, the problem admits an analytical solution (at least in the free-surface case, where the dynamics of the surrounding air is ignored). The solution has the form (see \cite{monaghan1994simulating})
    \begin{equation}
        \vv  = \begin{pmatrix}
            -\xi x\\
            \xi y
        \end{pmatrix},
        \quad 
        p = -\frac{\xi^2 \rho_\mathrm{water}}{\frac{1}{a^2}+ \frac{1}{b^2}}\left( \left(\frac{x}{a}\right)^2 + \left(\frac{y}{b}\right)^2 - 1\right), \quad \forall\xx \in \Omega_\mathrm{water},
    \end{equation}
    where 
    \begin{equation}
        \Omega_\mathrm{water} = \left \{\xx \in \Omega: \left(\frac{x}{a}\right)^2 + \left( \frac{y}{b} \right)^2 < 1\right \}
    \end{equation}
    and the triplet $a(t), b(t), \gamma(t)$ is a solution of an ordinary differential equation:
    \begin{equation}
        \begin{split}
            \dv{a}{t} &=  -\xi a,\\
            \dv{b}{t} &= \xi b,\\
            \dv{\xi}{t} &= \xi^2\frac{a^2 - b^2}{a^2 + b^2}
        \end{split}
    \end{equation}
    with initial condition $a(0) = b(0) = 1$ and $\xi(0) = \frac{1}{2}$. Comparing the numerical results with the semi-analytical solution allows us to perform a convergence analysis, as depicted in Figure \ref{fig:cpatch}. Despite using a two-phase model in place of a free surface condition, we see that the estimated order of convergence (EOS) is about 1 for both velocity field and pressure.
    
    \subsection{Dam break}
    
\begin{figure}[htb!]
    	\centering
    	\begin{subfigure}{0.33\linewidth}
    		\includegraphics[width=\linewidth]{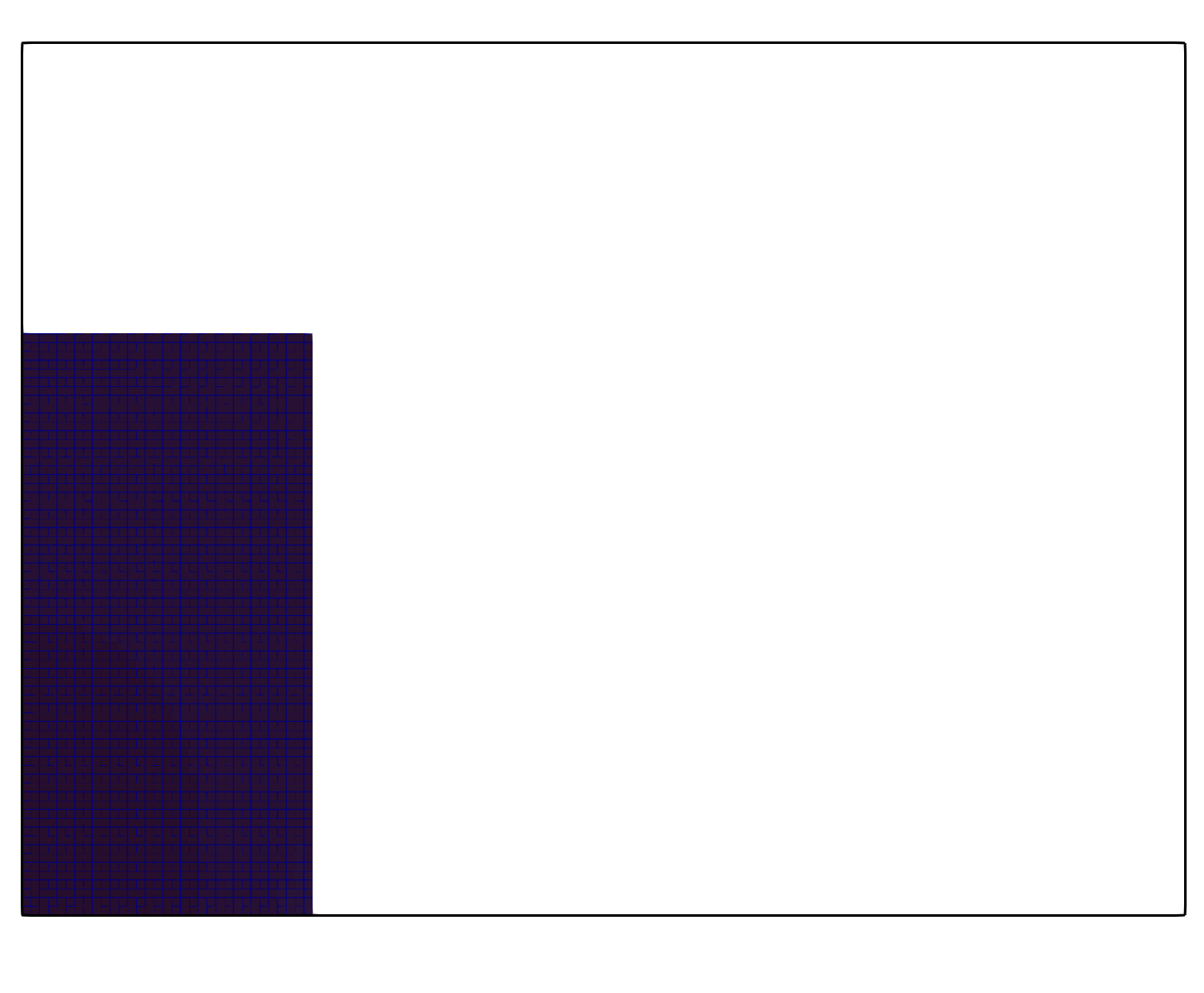}
    		\caption{$t = 0.0\unit{s}$}
    	\end{subfigure}%
    	\begin{subfigure}{0.33\linewidth}
    		\includegraphics[width=\linewidth]{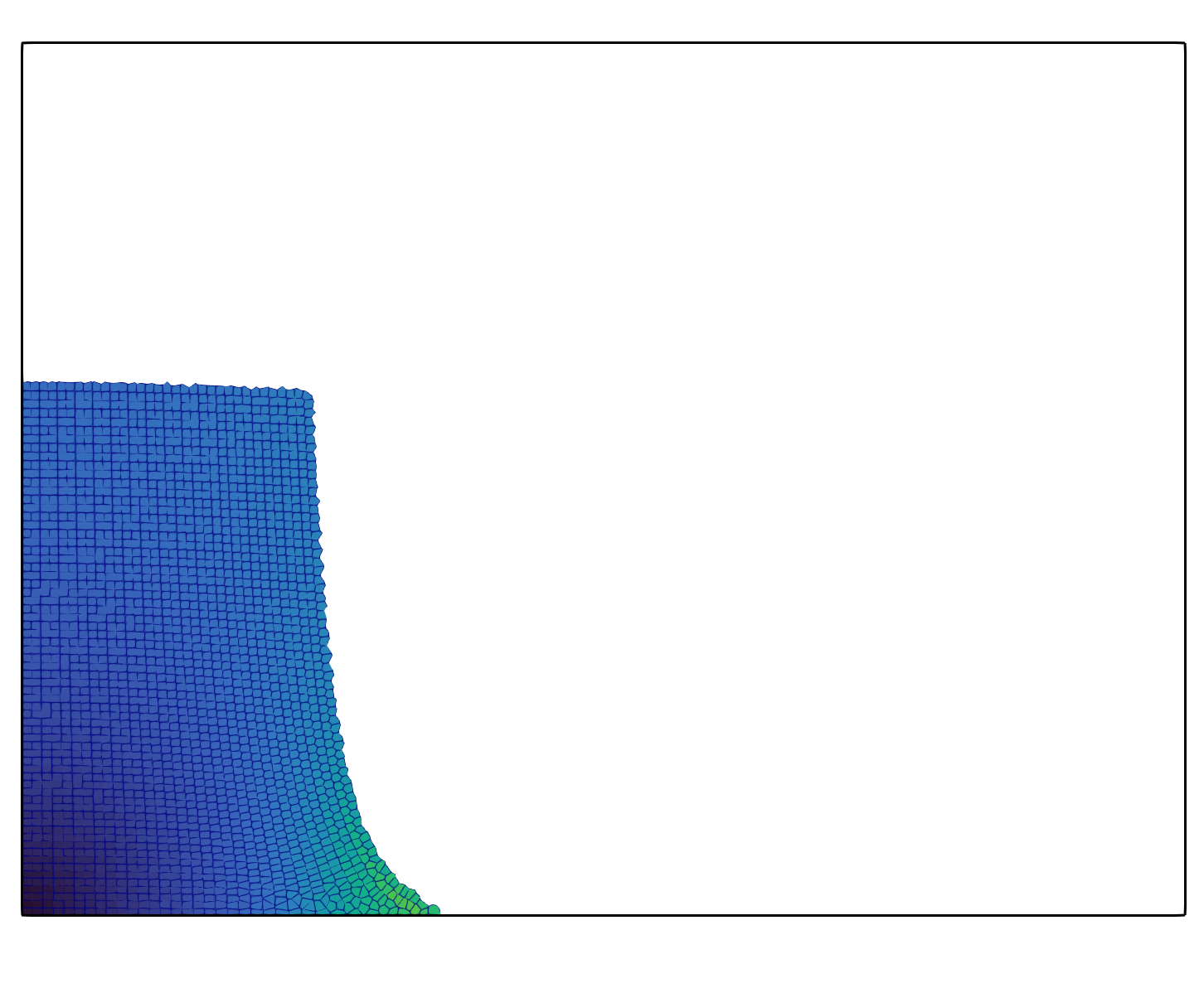}
    		\caption{$t = 0.2\unit{s}$}
    	\end{subfigure}%
    	\begin{subfigure}{0.33\linewidth}
    		\includegraphics[width=\linewidth]{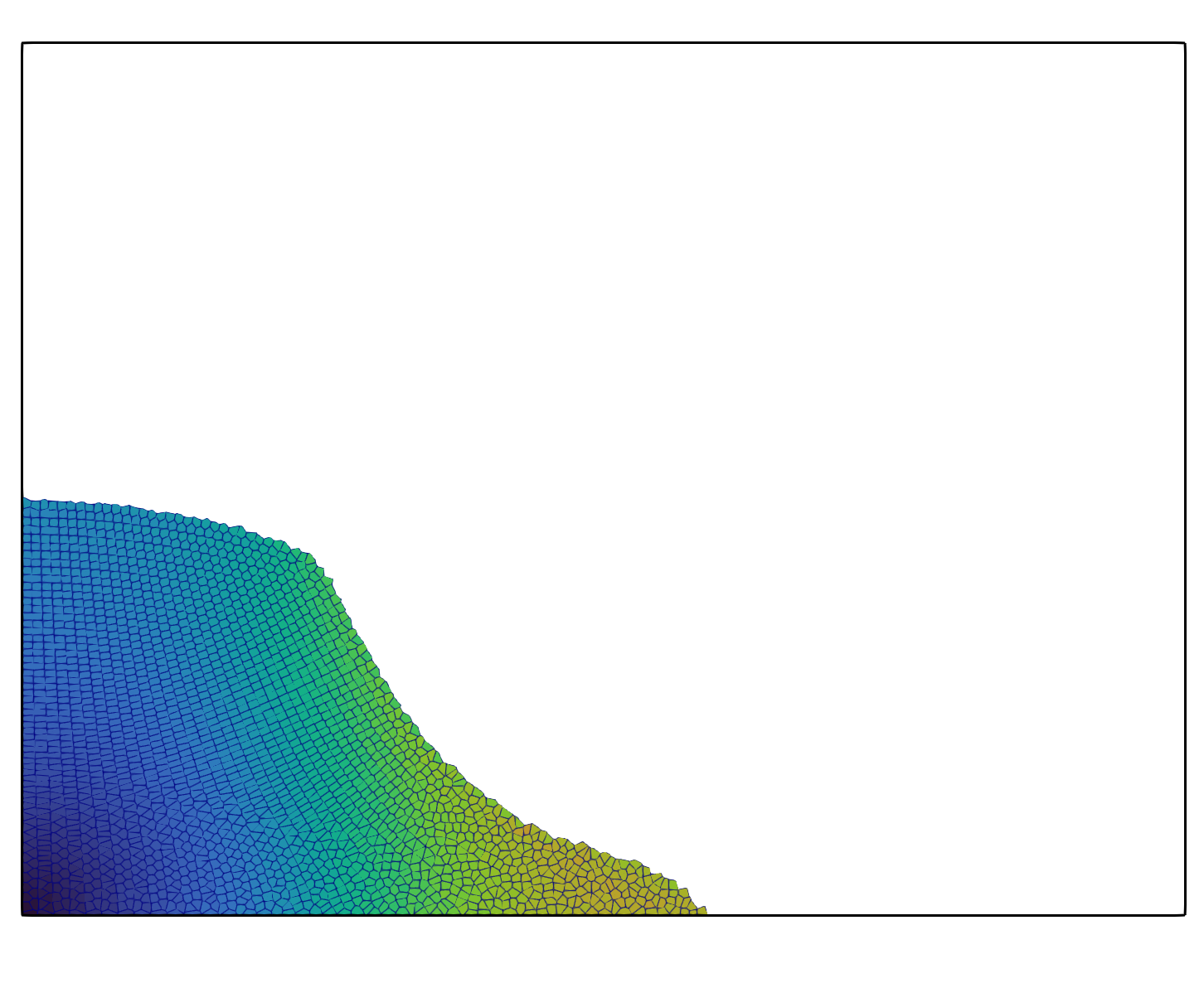}
    		\caption{$t = 0.4\unit{s}$}
    	\end{subfigure}
    	\begin{subfigure}{0.33\linewidth}
    		\includegraphics[width=\linewidth]{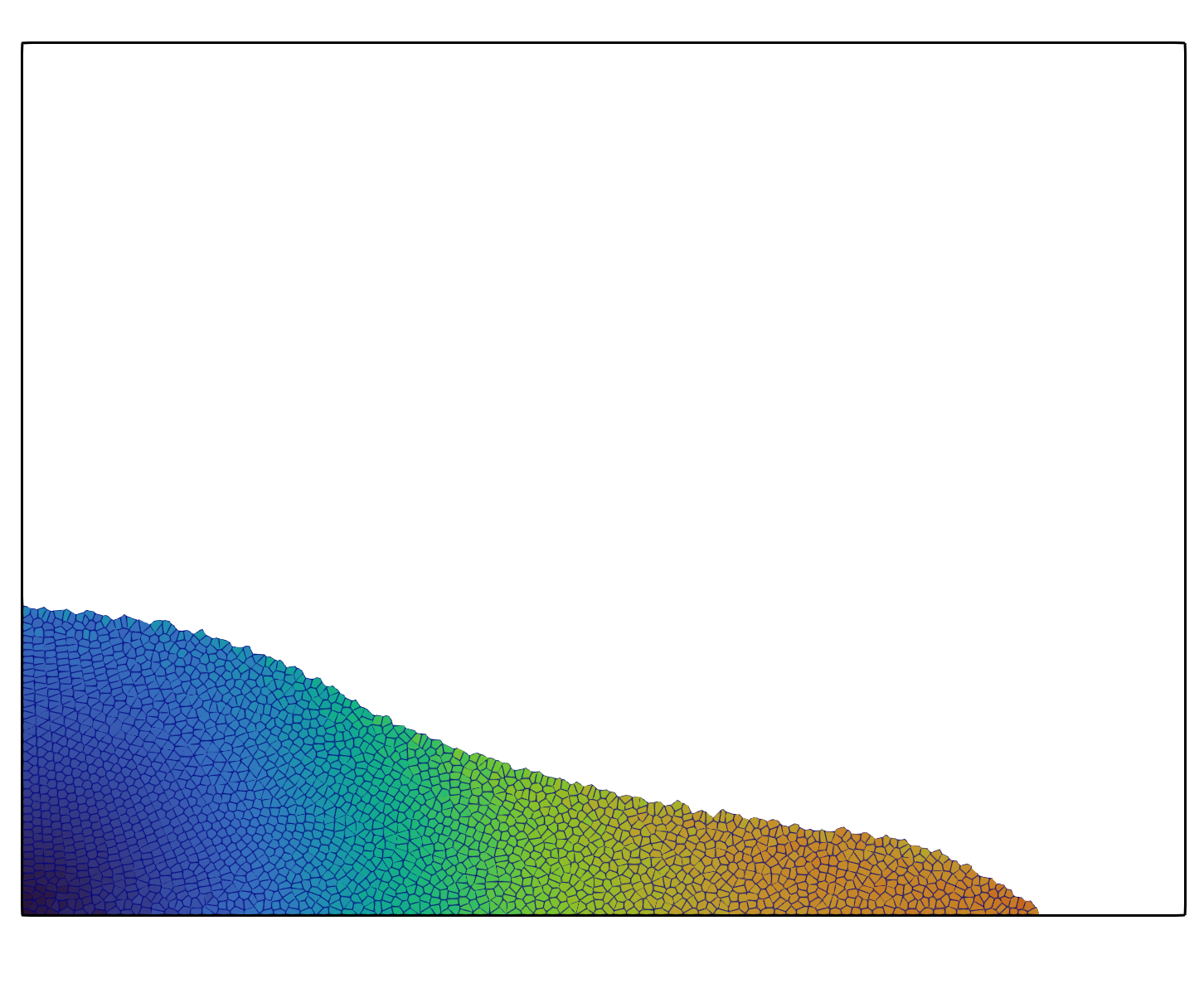}
    		\caption{$t = 0.6\unit{s}$}
    	\end{subfigure}%
    	\begin{subfigure}{0.33\linewidth}
    		\includegraphics[width=\linewidth]{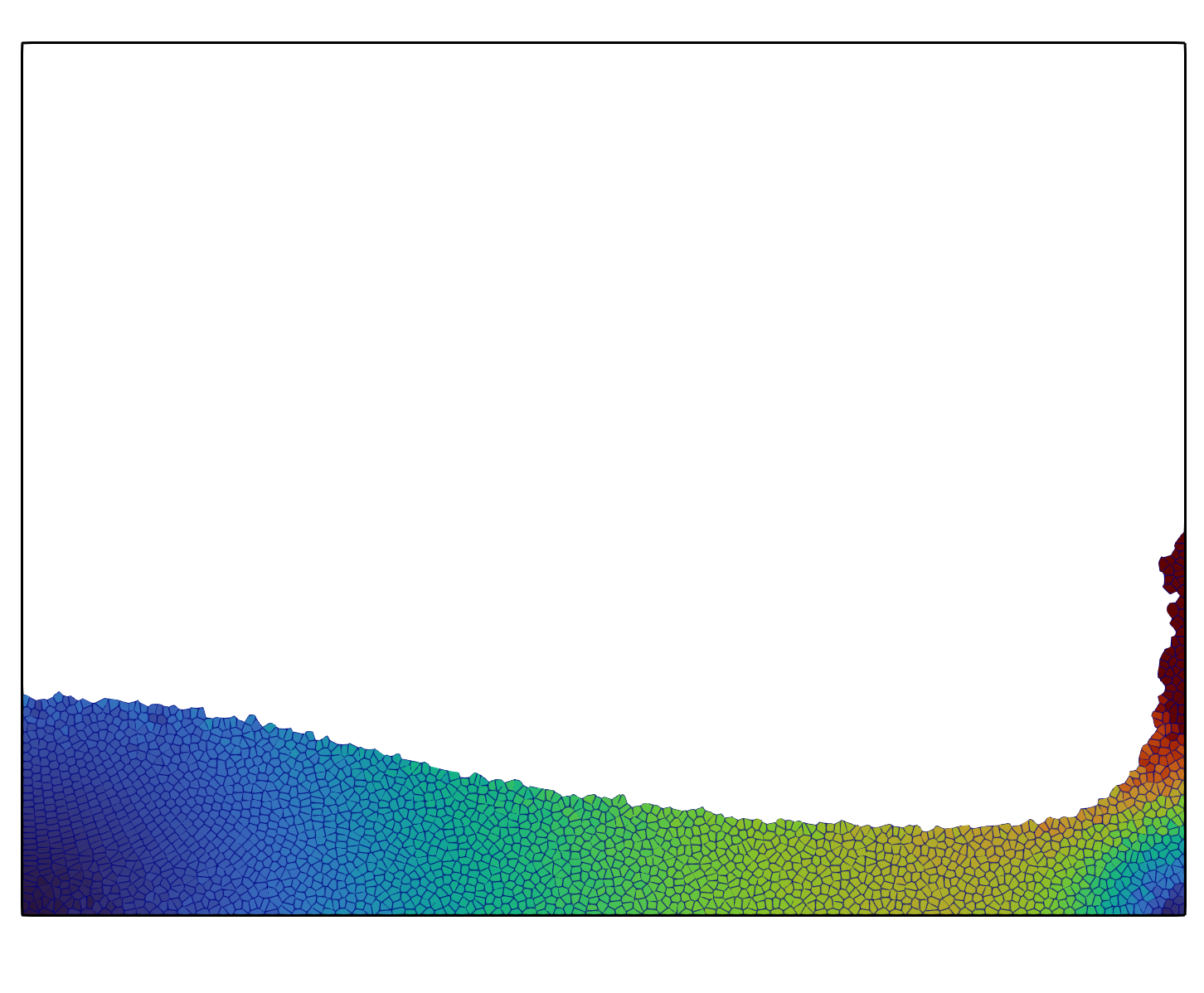}
    		\caption{$t = 0.8\unit{s}$}
    	\end{subfigure}%
    	\begin{subfigure}{0.33\linewidth}
    		\includegraphics[width=\linewidth]{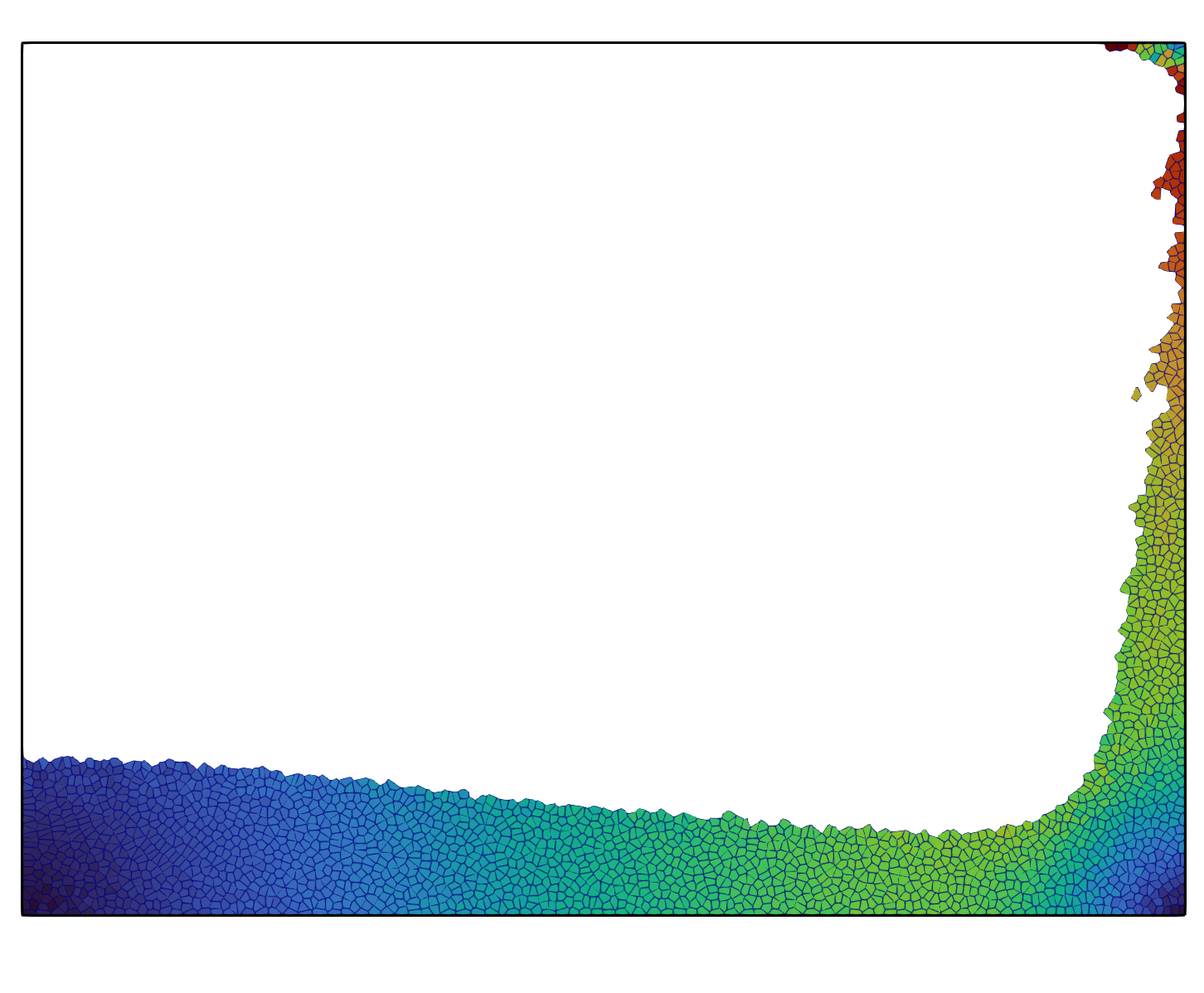}
    		\caption{$t = 1.0\unit{s}$}
    	\end{subfigure}
    	\centering
    	\includegraphics[width = 0.3\linewidth]{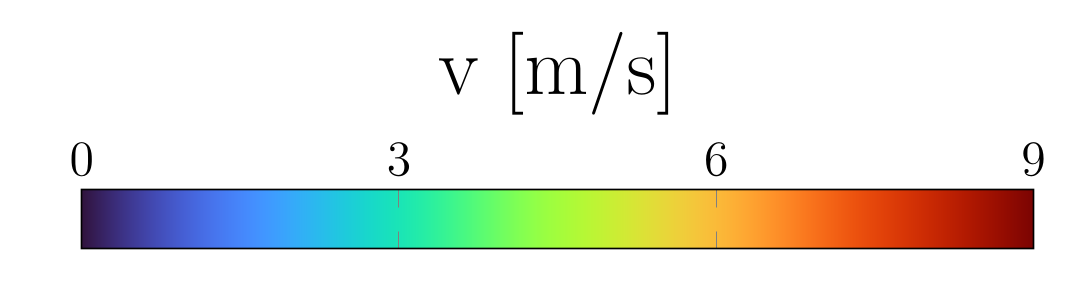}
    	\caption{Selected frames of the dam break test. The color map indicates magnitude of velocity. Air cells are hidden from view.}
    	\label{fig:dambreak_v}
    \end{figure}

\begin{figure}
        \centering
        \begin{subfigure}{0.5\linewidth}
            \includegraphics[width=\linewidth]{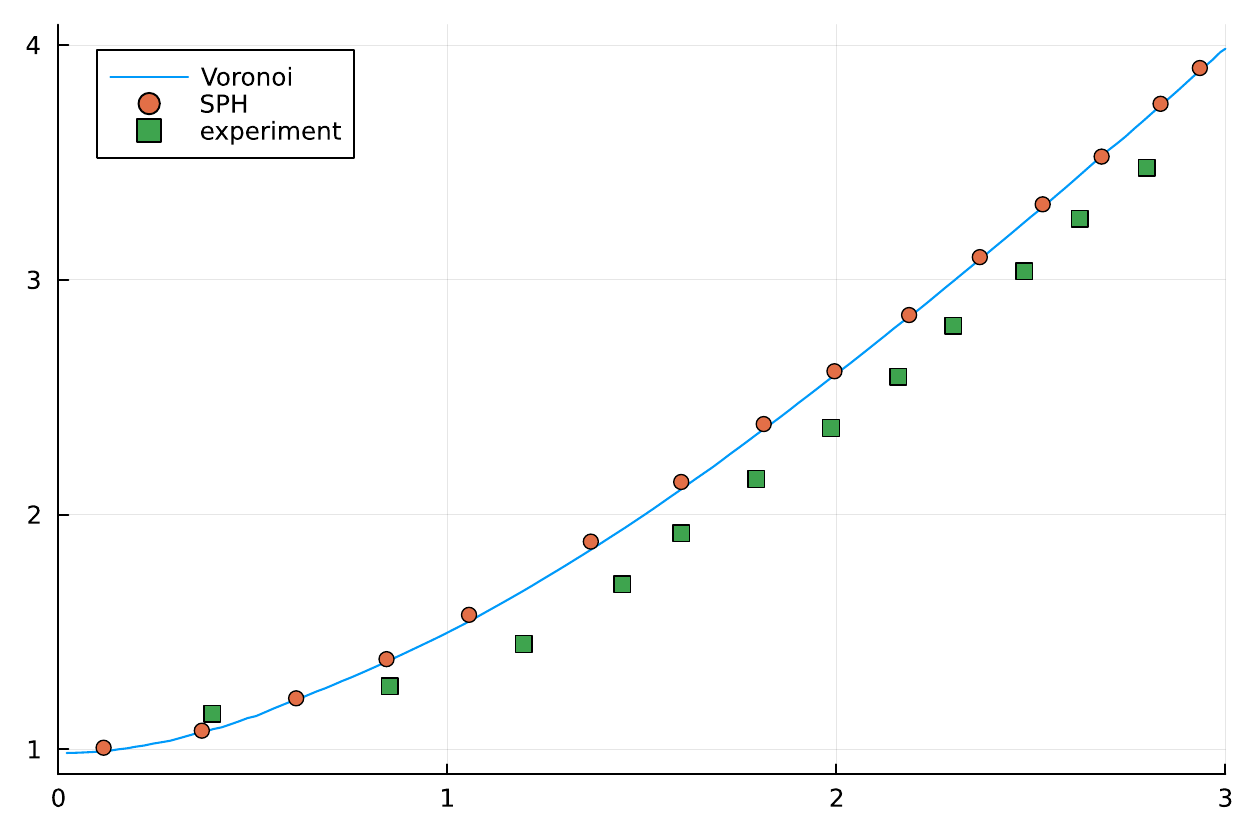}
            \caption{wavefront location $X$}
        \end{subfigure}%
        \begin{subfigure}{0.5\linewidth}
            \includegraphics[width=\linewidth]{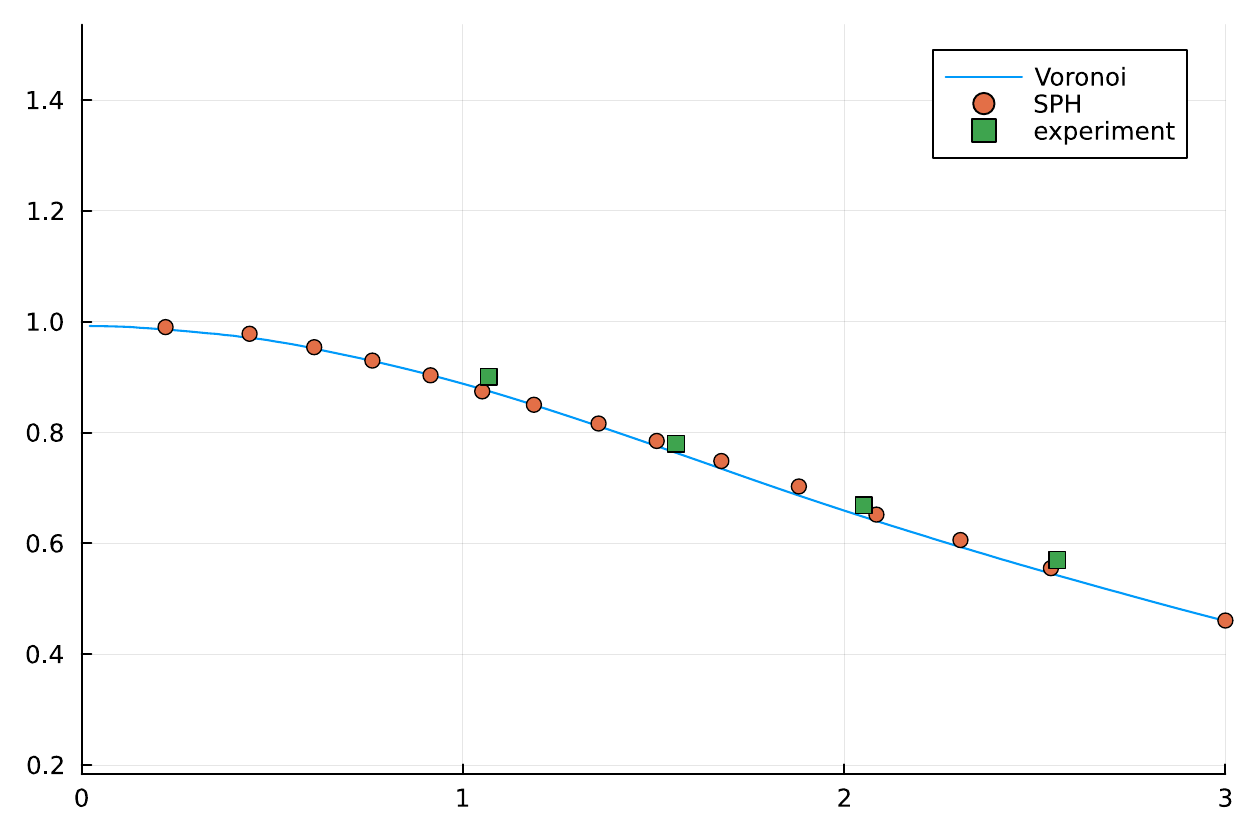}
            \caption{water column height $H$}
        \end{subfigure}
        \caption{Comparing our numerical simulation in dam break test with an SPH result \cite{violeau} and an experiment by Koshizuka and Oka \cite{koshizuka1996moving}. We observe excellent agreement with the particle simulation and small discrepancy with the experiment, pointing to the possibility that some physical aspect is missing in both simulations.}
        \label{fig:dambreak_graph}
    \end{figure}

    As for a more challenging test that features viscosity and gravity, we consider a two-dimensional dam break problem. The computational domain is
    $ \Omega = (0, 4\unit{m})\times (0,3\unit{m}) $.
    At $t=0$, a water column $\Omega_\mathrm{water} = (0,1\unit{m})\times (0,2\unit{m})$ is defined, which then loses the stability due to the vertical free surface, and generates a water wave that travels to the right until it collides against the right wall. The density of water and the surrounding air remain the same as in the circular patch test, but now, we also consider viscosity of water $\mu_\mathrm{water} = 8.9\cdot10^{-4} \unit{Pa.s}$, air viscosity $\mu_\mathrm{air} = 3.7\cdot 10^{-5} \unit{Pa.s}$, and a homogeneous gravitational field $g = 9.8 \unit{m/s^2}$. The no-slip boundary condition is prescribed at every wall. We monitor the evolution of the following two quantities:
    \begin{align}
        X = \sup_{\xx\in \Omega_\mathrm{water}}\frac{x }{a}, \quad H = \sup_{y\in \Omega_\mathrm{water}}\frac{y}{b},
    \end{align}
    that represent the dimensionless $x$-coordinate of the wavefront and the dimensionless height of the water column over a dimensionless time
    \begin{equation}
        T = \sqrt{\frac{2gt^2}{b}},
    \end{equation}
    where $a = 1\unit{m}$ and $b = 2\unit{m}$ are the initial width and height of the water column, respectively. Results up to $t = 1\unit{s}$ are depicted in the Figures \ref{fig:dambreak_v} and \ref{fig:dambreak_graph}. We observe a good agreement with an experimental measurement \cite{koshizuka1996moving} and an excellent agreement with an SPH reference solution.
    
    \subsection{Shock-bubble interaction}
\begin{figure}[htb!]
        \centering
        \includegraphics[width=0.3\linewidth]{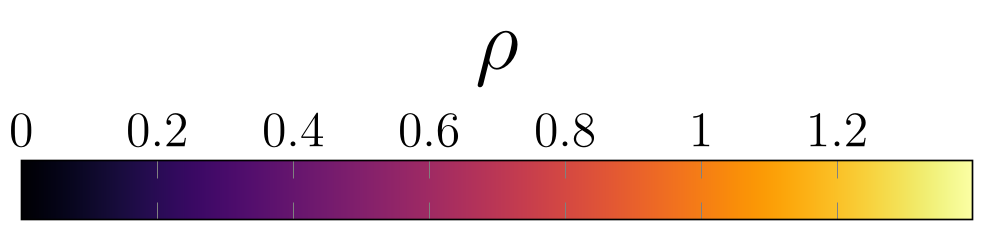}
        \begin{subfigure}{0.5\linewidth}
            \includegraphics[width=\linewidth]{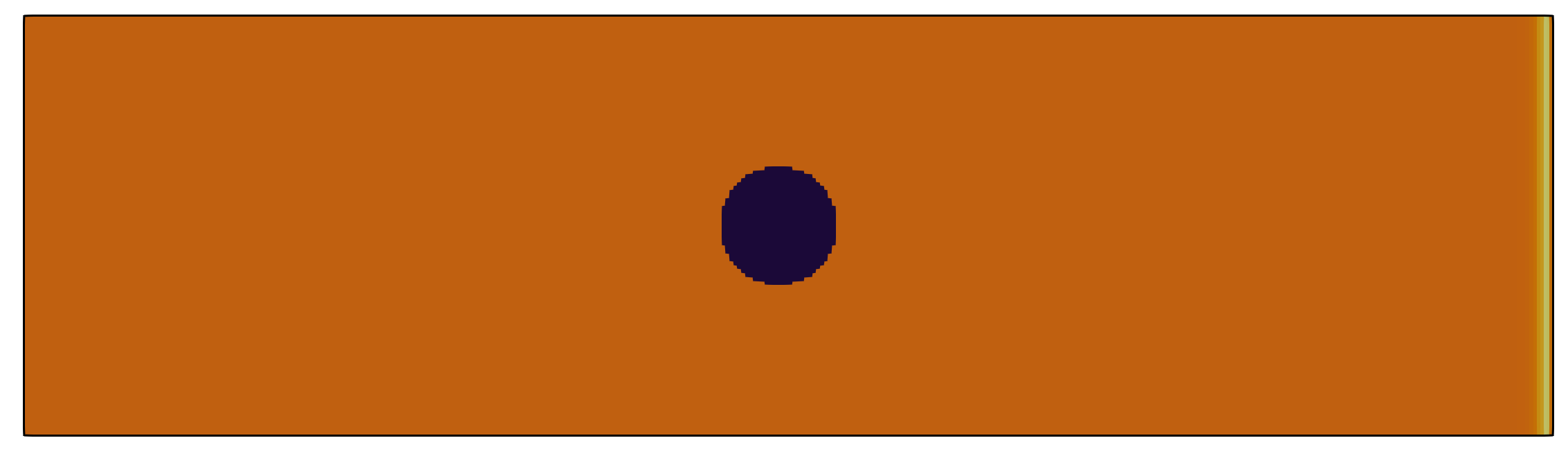}
            \caption{$t = 0$}
        \end{subfigure}%
        \begin{subfigure}{0.5\linewidth}
            \includegraphics[width=\linewidth]{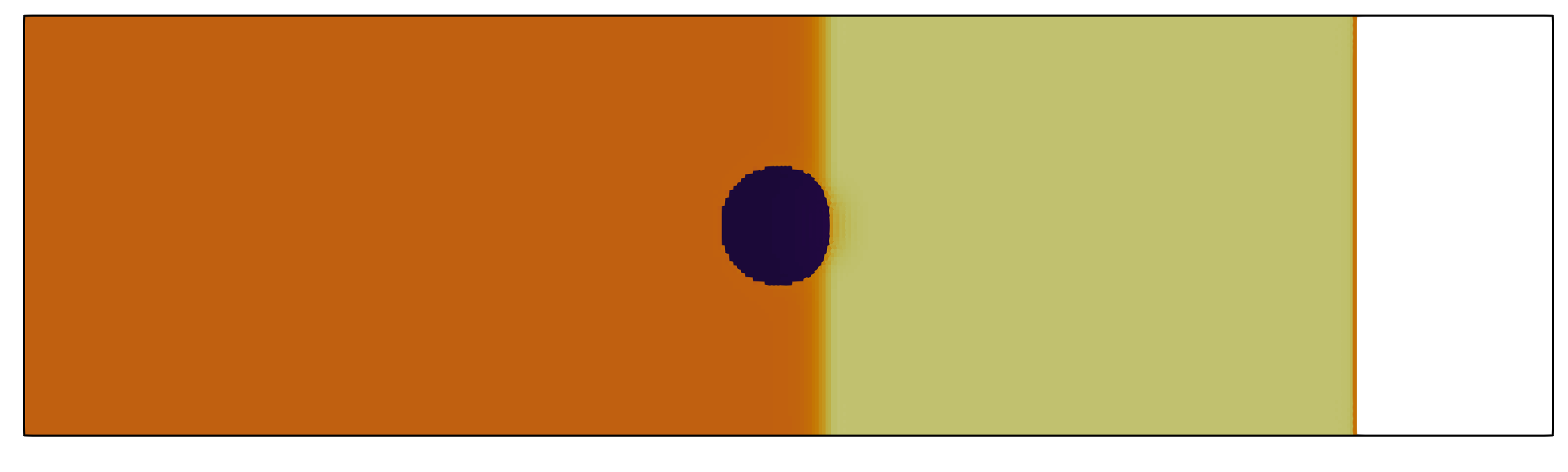}
            \caption{$t = \frac{1}{2}t_\mathrm{end}$}
        \end{subfigure}
        \begin{subfigure}{0.5\linewidth}
            \includegraphics[width=\linewidth]{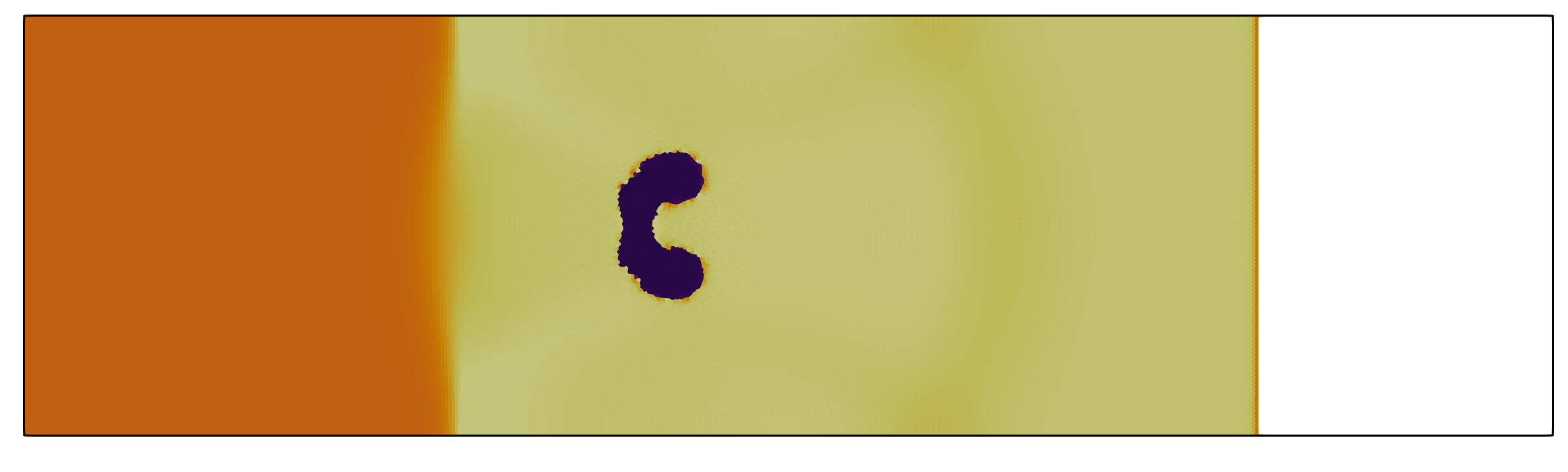}
            \caption{$t = \frac{3}{4}t_\mathrm{end}$}
        \end{subfigure}%
        \begin{subfigure}{0.5\linewidth}
            \includegraphics[width=\linewidth]{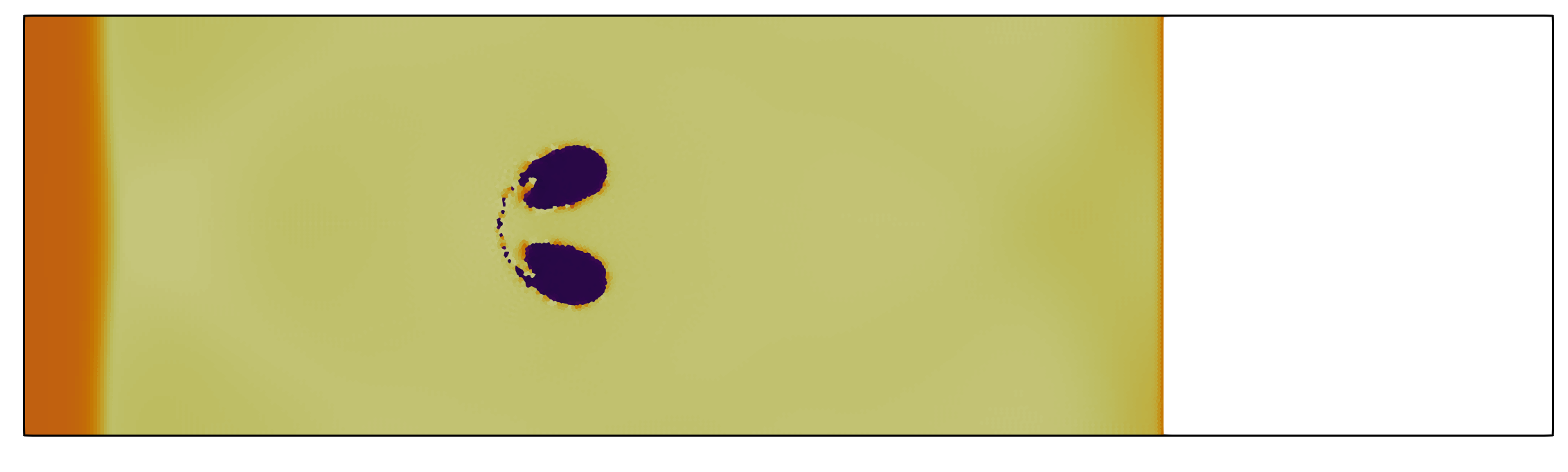}
            \caption{$t = t_\mathrm{end}$}
        \end{subfigure}
        \begin{subfigure}{\linewidth}
                \centering
             \includegraphics[width=\linewidth]{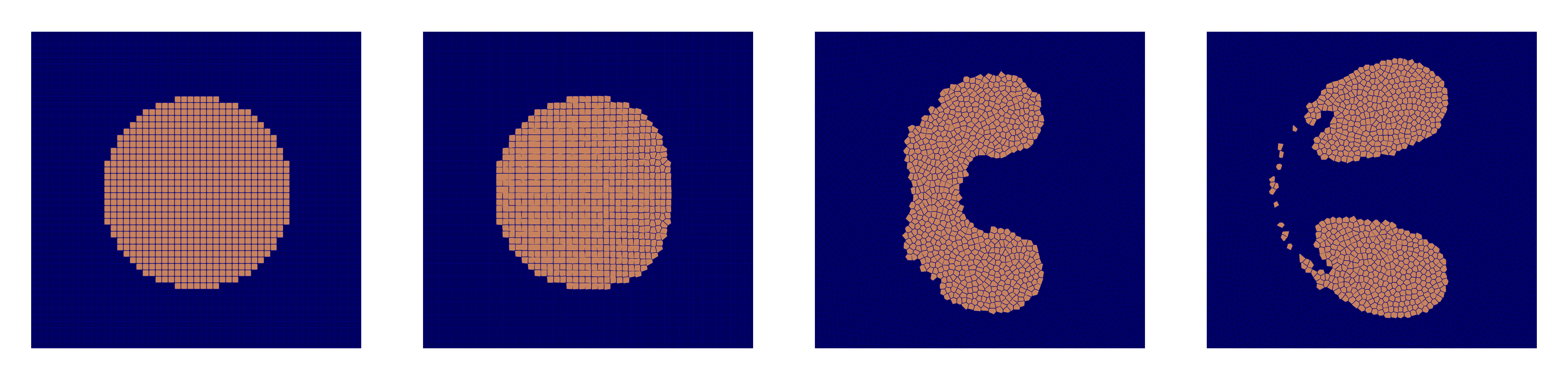}
            \caption{Shape of the bubble in detail at $t=0$, $\frac12 t_\mathrm{end}$, $\frac34 t_\mathrm{end}$, and $t_\mathrm{end}$.}
        \end{subfigure}
        \caption{The shock-bubble interaction at selected time frames. Colormap indicates density.}
        \label{fig:shock-bubble}
    \end{figure}

    As a compressible test, a two-dimensional shock-bubble interaction benchmark is presented. In the initial configuration, a circular bubble of helium 
        \begin{equation}
            \Omega_\mathrm{He} = \{ \xx \in \Omega : (x-0.32)^2 + (y-0.089)^2 < 0.025^2\}
        \end{equation}
    with density $\rho_\mathrm{He} = 0.182 \unit{kg/m^3}$ and adiabatic index $\gamma_\mathrm{He} = 1.648$ is surrounded by air in a piston
    \begin{equation}
        \Omega = (0,0.65)\times (0,0.178) .
    \end{equation}
    The density of air is $\rho_\mathrm{air}=1 \unit{Kg/m^3}$ and the adiabatic index is $\gamma = 1.4$. Both fluids, modeled as ideal fluids, are initially at rest with pressure $P = 10^5 \unit{Pa}$. The right hand side of the piston is moving inward with velocity $v = 124.824 \unit{m/s}$. This generates a shock, which deforms the helium bubble. To account for entropy production near the shock wave, we add Stone-Norman artificial viscosity \cite{stone1992zeus}:
    	\begin{equation}
		\mu^\text{art}_i =  \begin{cases}
			-(\delta r)^2 \rho_i (\Tr \mathbb{D}_i), & \Tr \mathbb{D}_i < 0\\
			0, & \Tr \mathbb{D}_i \geq 0
		\end{cases}.
		\label{eq:artificial_visc}
	\end{equation}	
    The results of the simulation, depicted in Figure \ref{fig:shock-bubble}, are for qualitative comparison with the reference \cite{loubere2010reale}. 

    \subsection{Rotating square}
    \begin{figure}[htb!]
        \centering
        \begin{subfigure}{0.33\linewidth}
            \includegraphics[width=\linewidth]{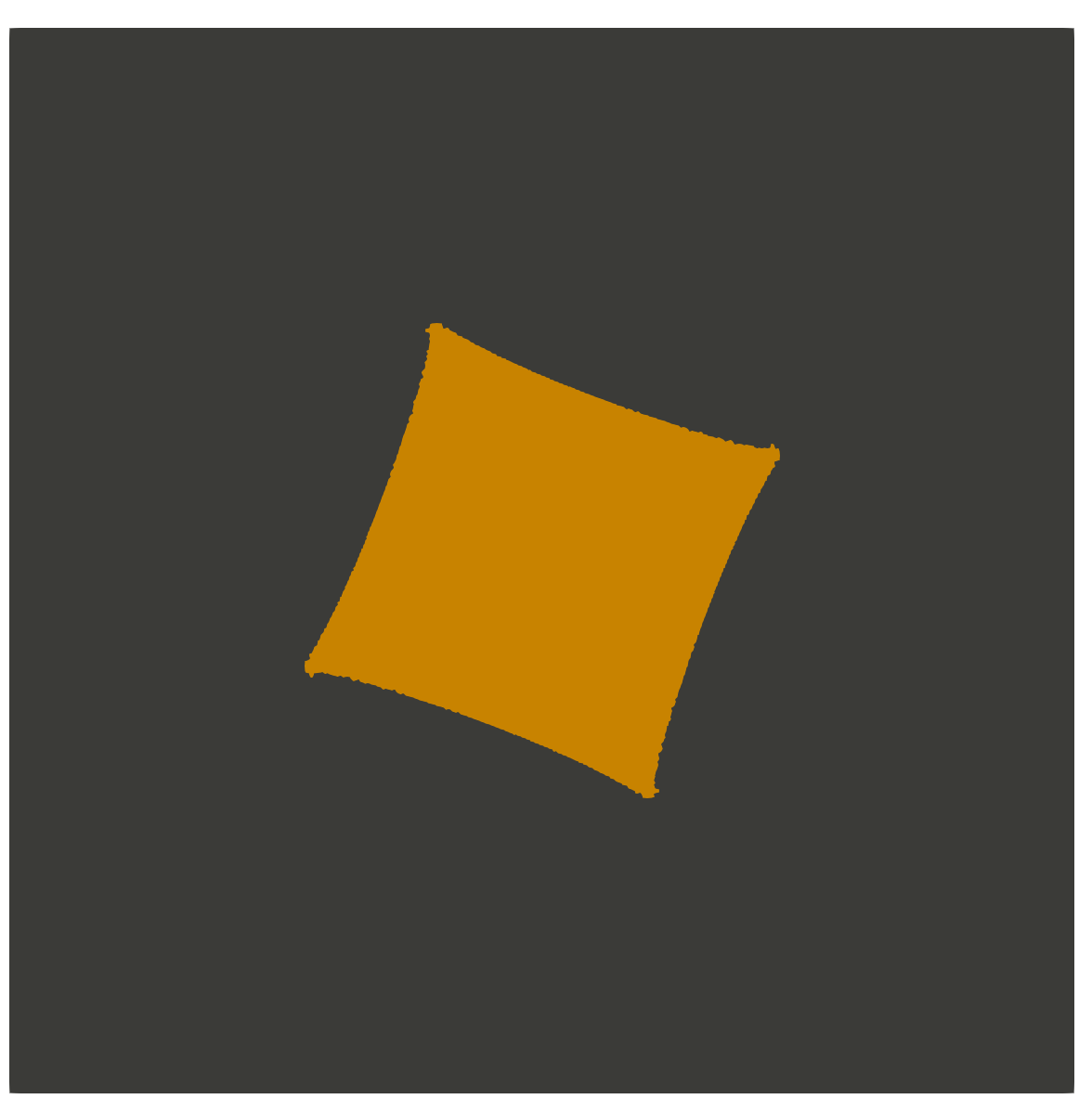}
            \caption{$t = 0.4$}
        \end{subfigure}%
        \begin{subfigure}{0.33\linewidth}
            \includegraphics[width=\linewidth]{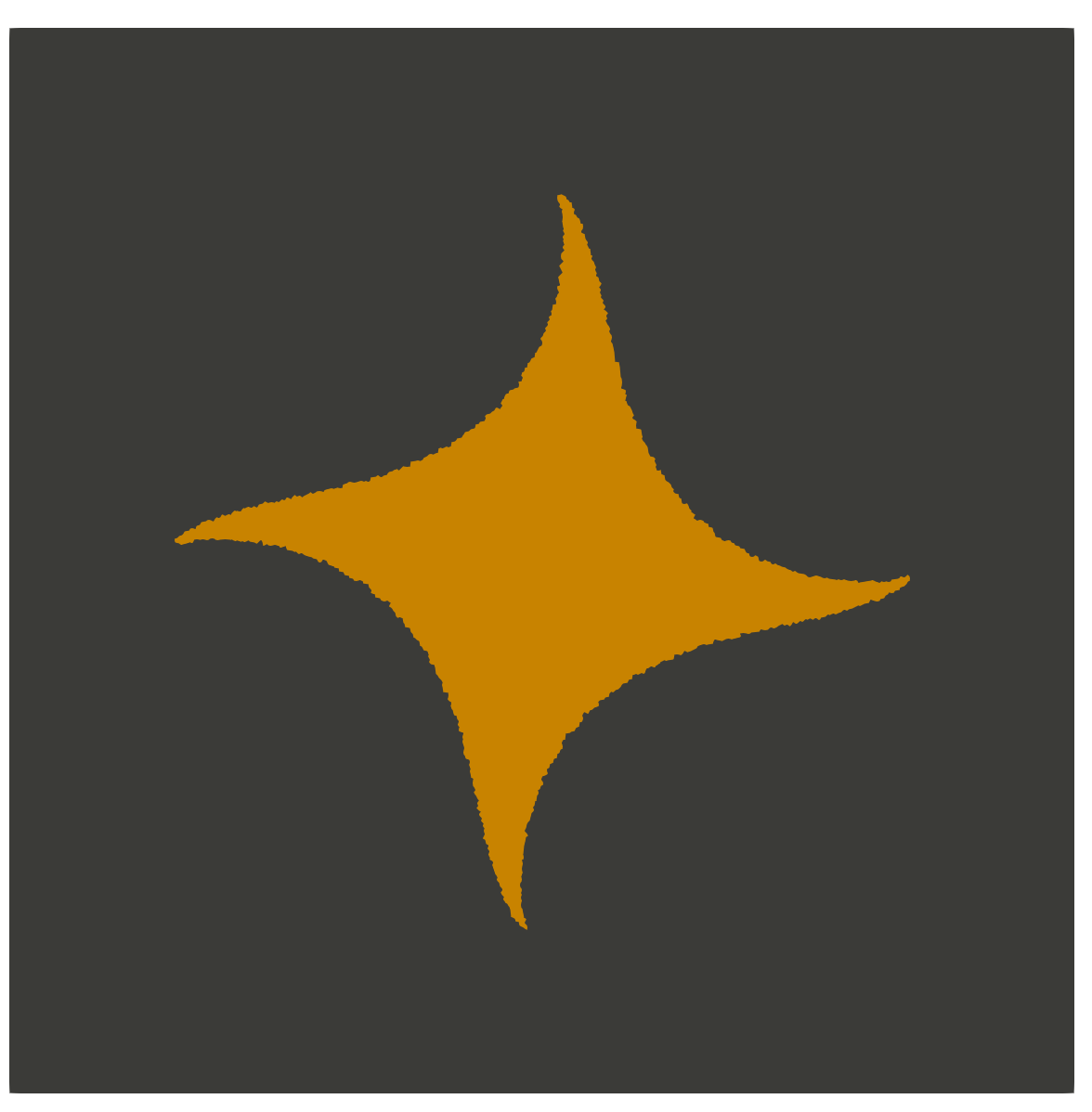}
            \caption{$t= 1.2$}
        \end{subfigure}%
        \begin{subfigure}{0.33\linewidth}
            \includegraphics[width=\linewidth]{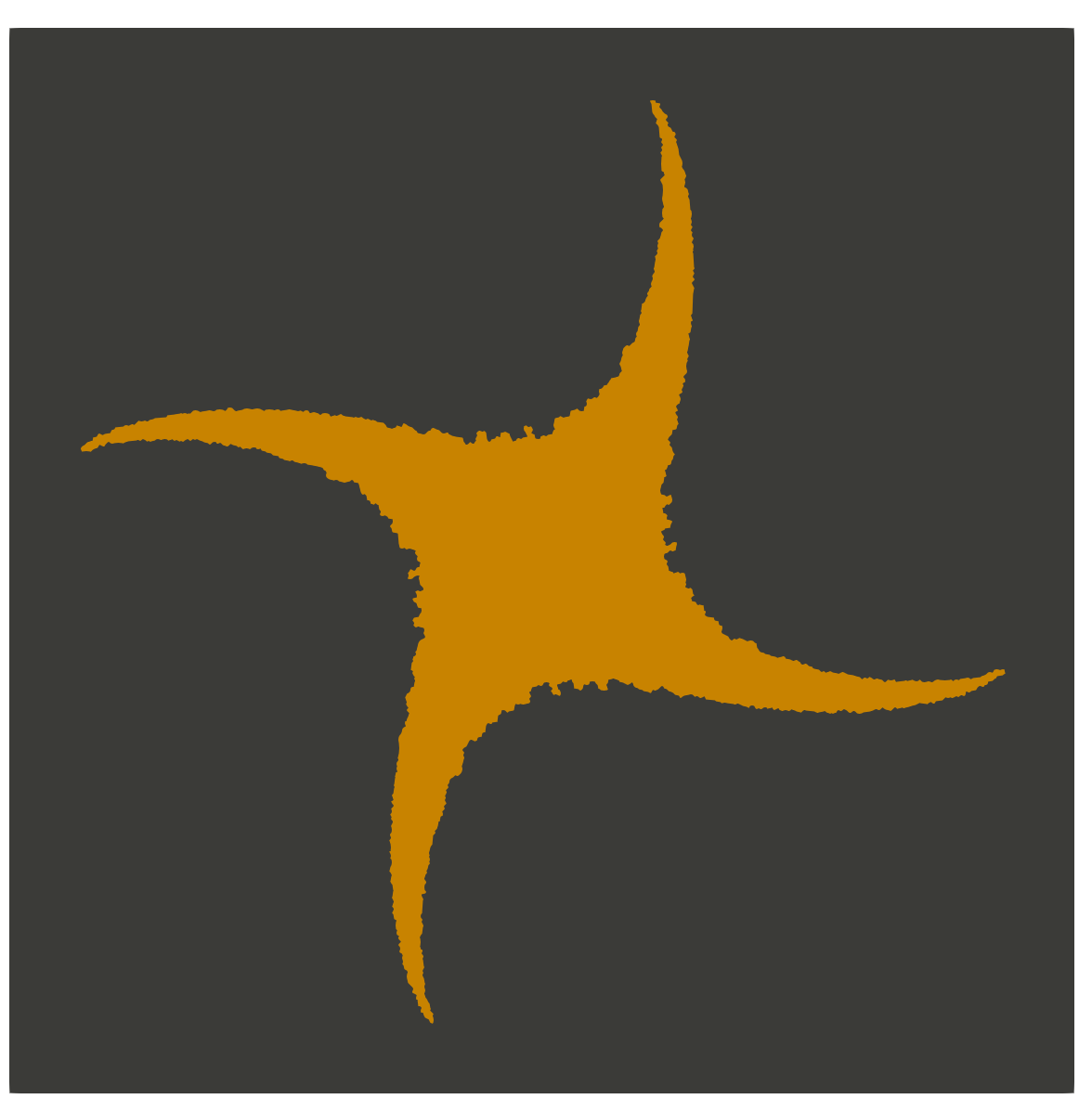}
            \caption{$t = 1.8$}
        \end{subfigure}
        \caption{Phase distribution.}
        \label{fig:rotsquare}
    \end{figure}
\begin{figure}[htb!]
        \centering
        \begin{subfigure}{0.33\linewidth}
            \includegraphics[width=\linewidth]{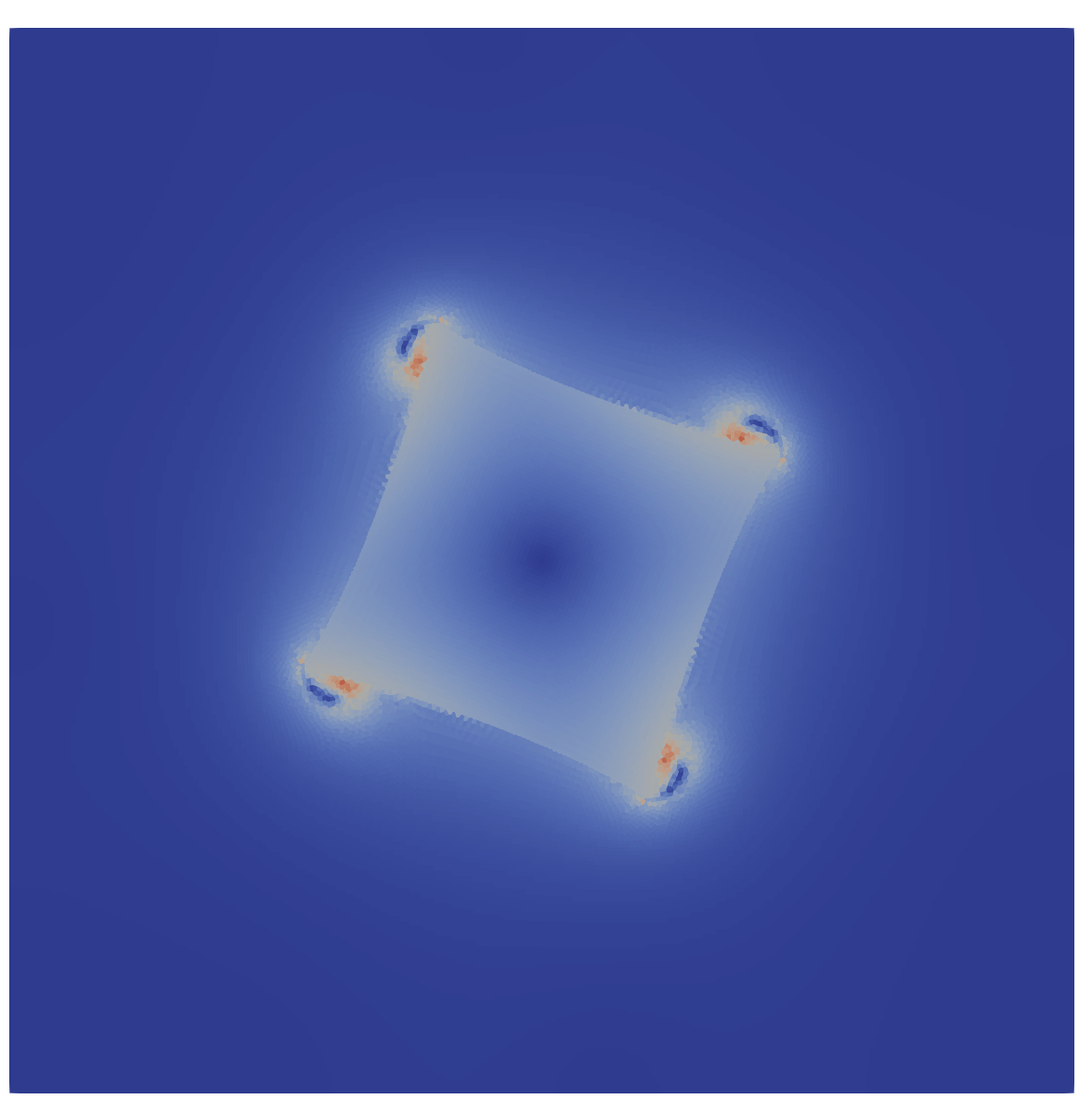}
            \caption{$t = 0.4$}
        \end{subfigure}%
        \begin{subfigure}{0.33\linewidth}
            \includegraphics[width=\linewidth]{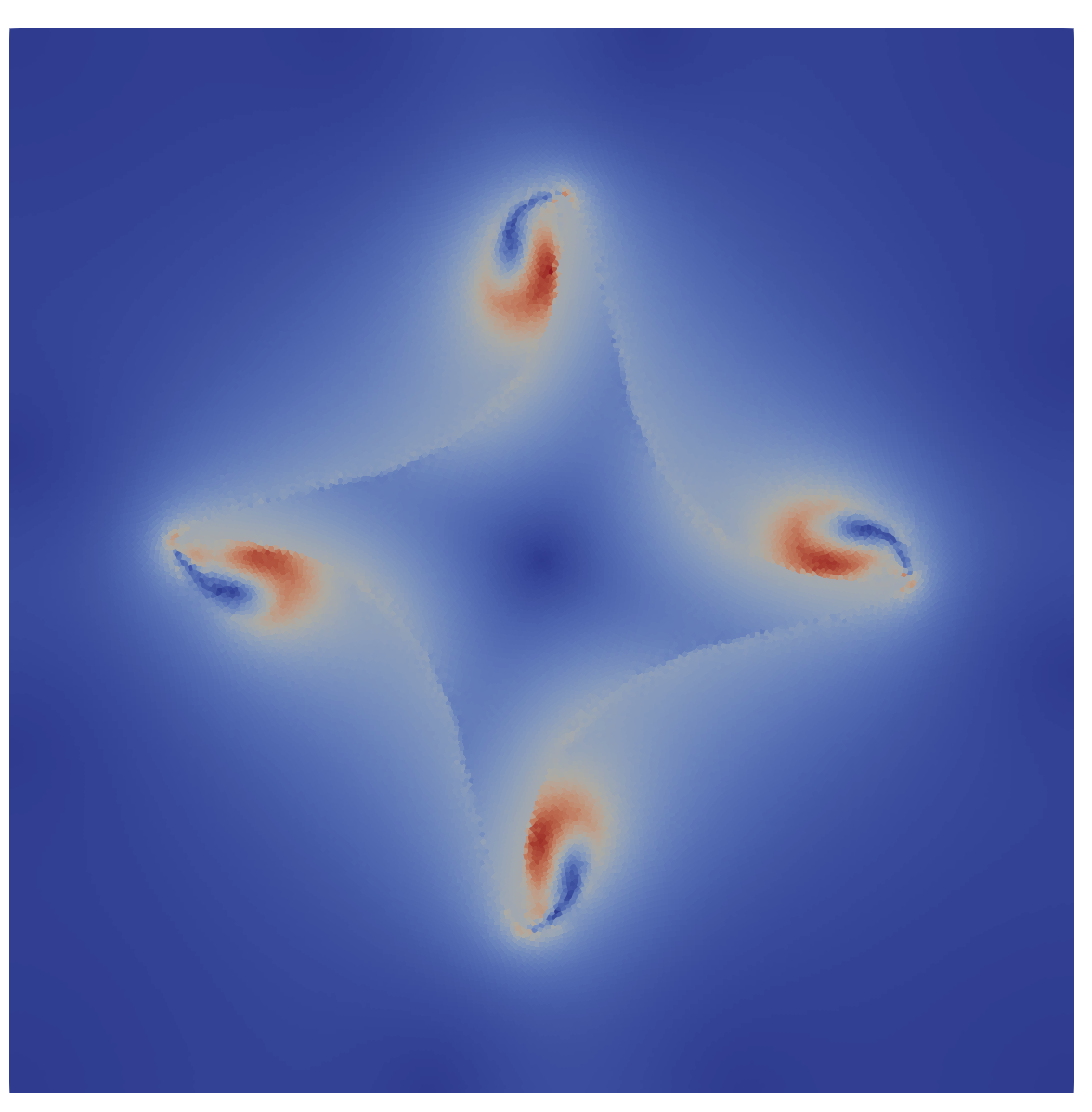}
            \caption{$t= 1.2$}
        \end{subfigure}%
        \begin{subfigure}{0.33\linewidth}
            \includegraphics[width=\linewidth]{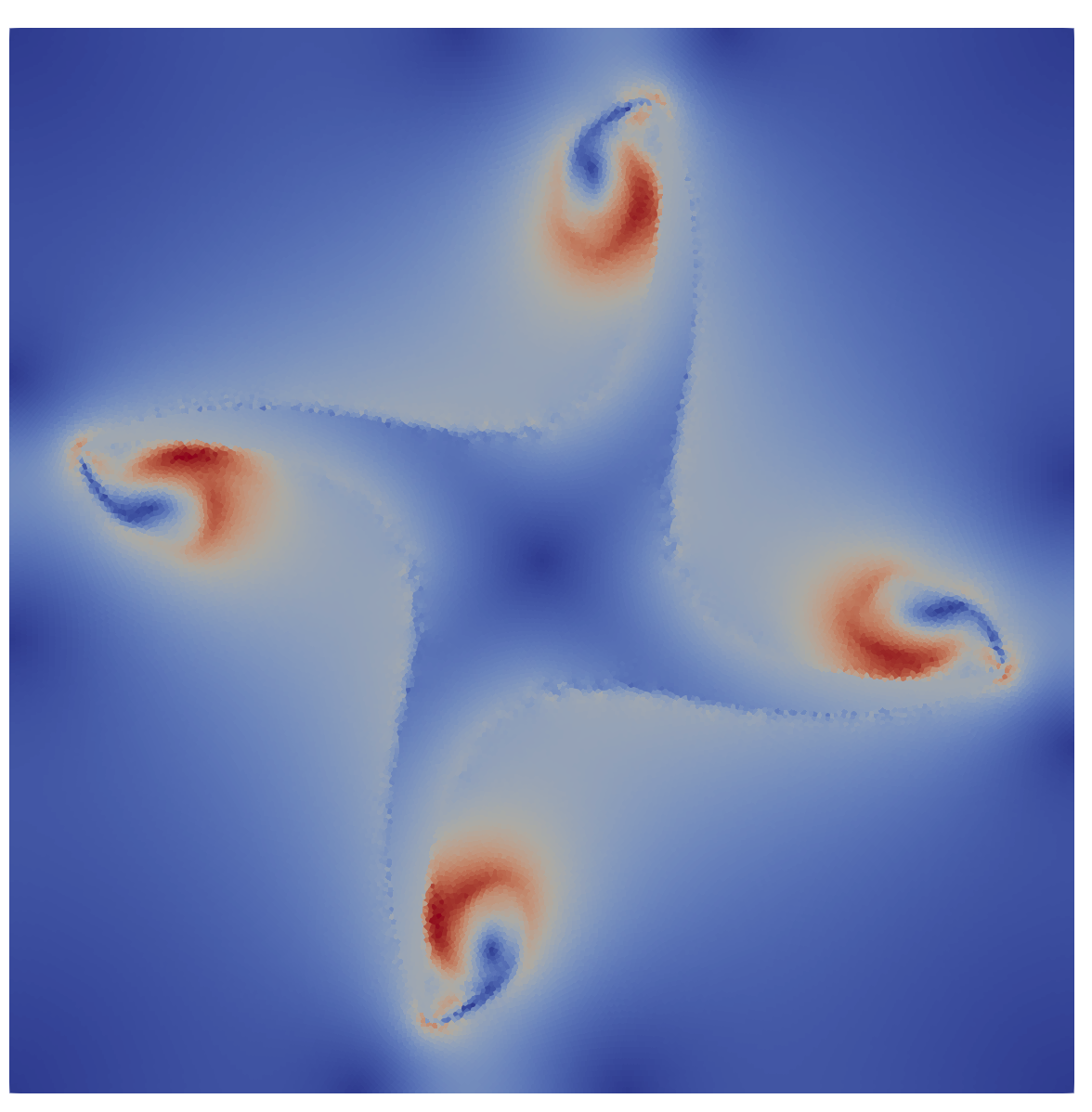}
            \caption{$t = 1.8$}
        \end{subfigure}
        \includegraphics[width=0.3\linewidth]{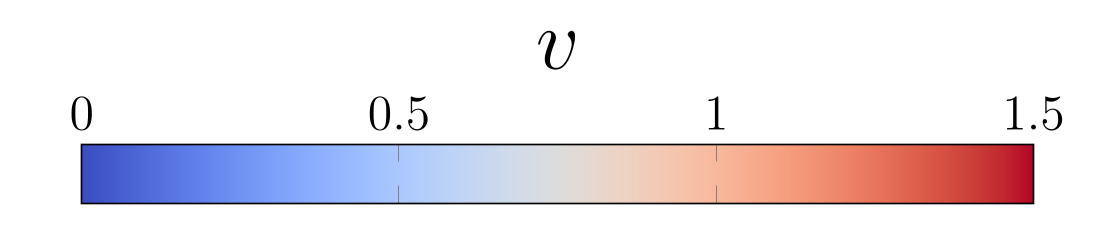}
        \caption{Magnitude of velocity.}
        \label{fig:rotsquare_v}
    \end{figure}

    In this free surface test, an incompressible liquid square
    \begin{equation}
        \Omega_\mathrm{water} = \left( -0.5,0.5\right)^2
    \end{equation}
    is subjected to a rigid body rotation:
    \begin{equation}
        \left.\vv \right|_{t=0} = \begin{pmatrix}
            -y \\
            x
        \end{pmatrix}. 
    \end{equation}
    Since the square is not rigid, the velocity field evolves and deforms the body into a star-like shape. This test was originally conceived as a free surface problem, but we shall simulate the surrounding air, with physical parameters identical to the circular patch benchmark. The size of the domain should be sufficiently large to accommodate the expanding liquid body. Our results, depicted in Figures \ref{fig:rotsquare} and \ref{fig:rotsquare_v}, are in qualitative agreement with a particle simulation \cite{sun2019consistent}. Unlike in SPH, we are not troubled by tensile instability.

    \subsection{Rising bubble}
    To validate the discretization of surface tension phenomena, let us study the classical rising bubble test. Let $\Omega = (0,1)\times(0,2)$ be the computational domain. As before, we consider a two-phase flow between water and air, but their roles are reversed in this case, with a small air body 
    \begin{equation}
        \Omega_\mathrm{air} = \left\{\xx \in \Omega: (x-0.5)^2 + (y-0.5)^2 < 0.25^2 \right\}
    \end{equation}
    surrounded by water. Free-slip boundary condition is prescribed at the lateral sides of the domain and a no-slip boundary condition at the top and bottom.
    
\begin{figure}[htb!]
        \centering
        \begin{subfigure}{0.5\linewidth}
            \includegraphics[width=\linewidth]{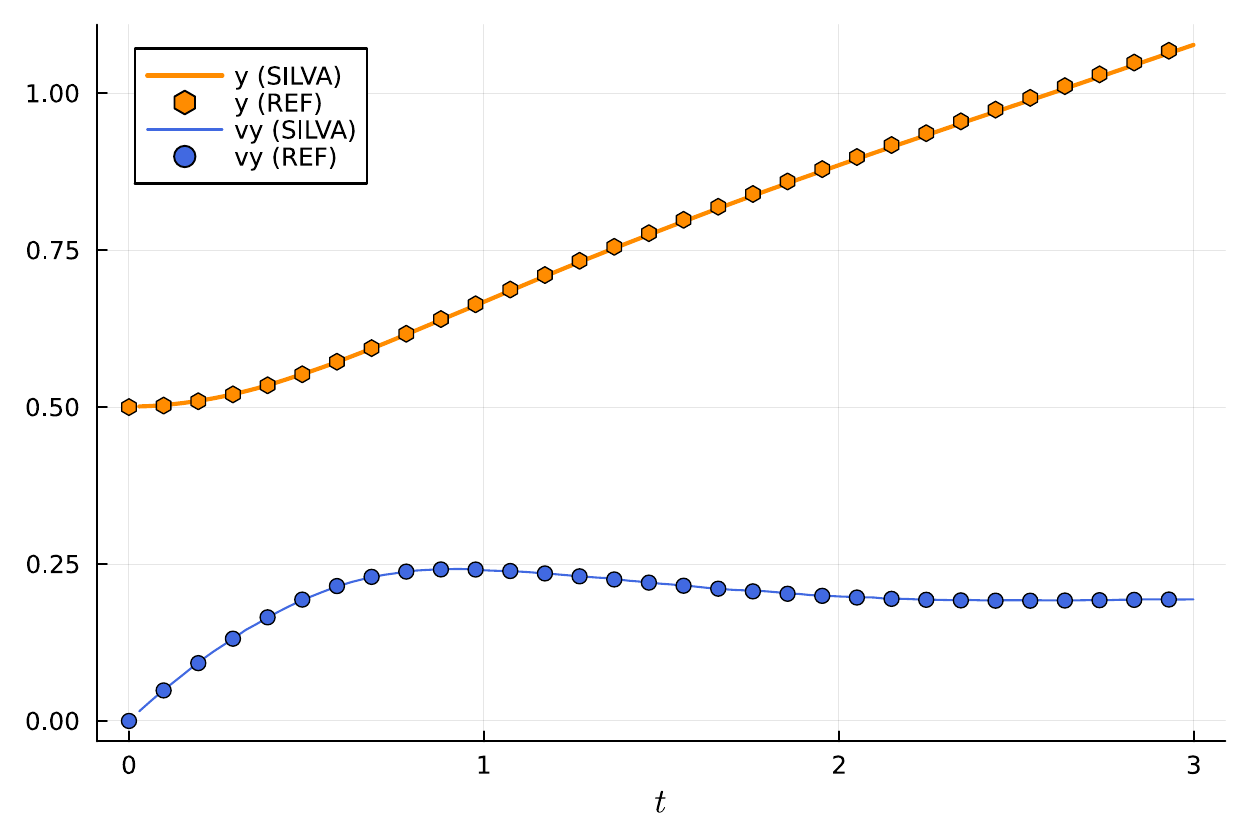}
            \caption{Vertical coordinate (y) and rising speed (vy) of the bubble.}
        \end{subfigure}%
        \begin{subfigure}{0.5\linewidth}
            \includegraphics[width=\linewidth]{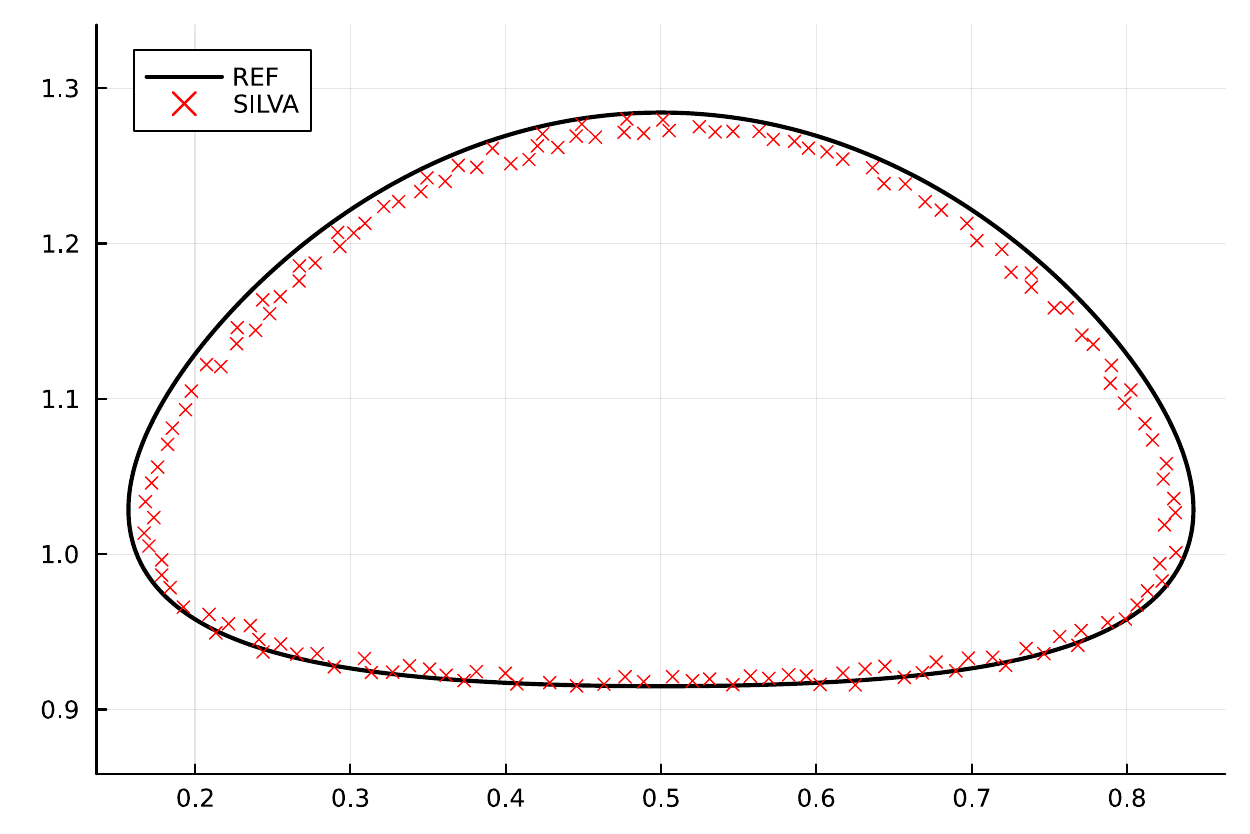}
            \caption{Bubble shape at $t= 3$.}
        \end{subfigure}
        \caption{Rising bubble (first setup). Comparison of our method (SILVA) with a reference solution (REF) based on finite element method \cite{hysing2009quantitative}. The resolution is 20 nodes per initial radius of the bubble.}
        \label{fig:bubble1}
    \end{figure}

\begin{figure}[htb!]
        \centering
        \begin{subfigure}{0.5\linewidth}
            \includegraphics[width=\linewidth]{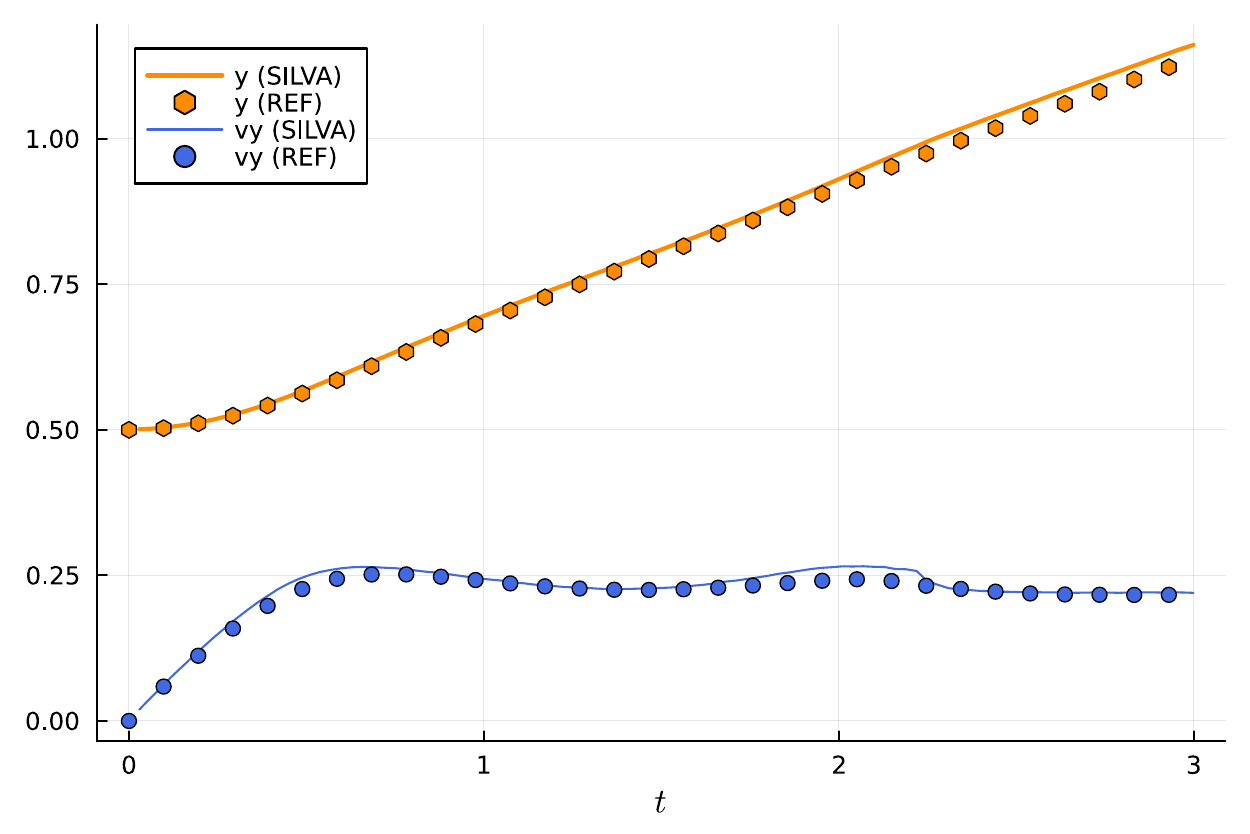}
            \caption{Vertical coordinate (y) and rising speed (vy) of the bubble.}
        \end{subfigure}%
        \begin{subfigure}{0.5\linewidth}
            \includegraphics[width=\linewidth]{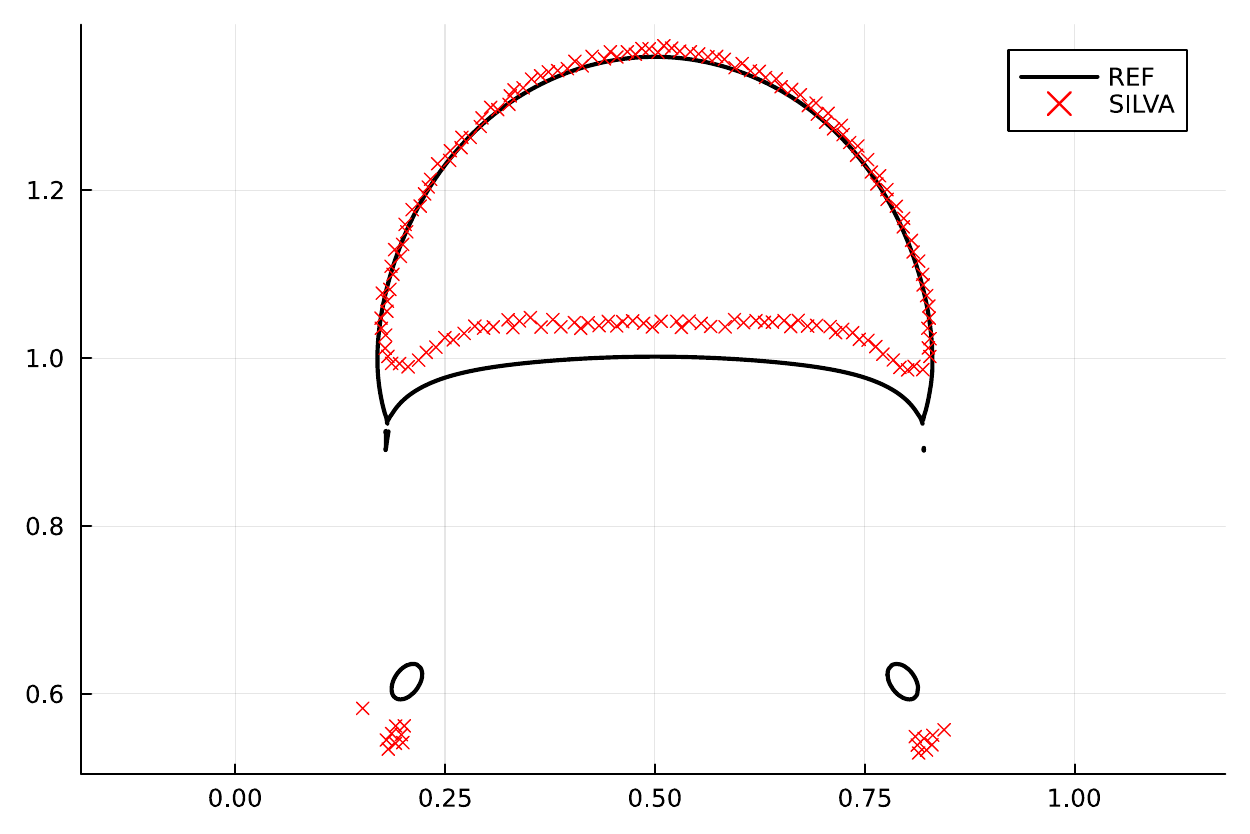}
            \caption{Bubble shape at $t= 3$.}
        \end{subfigure}
        \caption{Rising bubble (second setup).}
        \label{fig:bubble2}
    \end{figure}

    In the first setup, the physical parameters are: $\rho_\mathrm{water} = 1000$, $\rho_\mathrm{air} = 100$, $\mu_\mathrm{water} = 10$, $\mu_\mathrm{air} = 1$, $\sigma = 24.5$ and $g = 0.98$. Two properties of the bubble are measured over time: the $y$-coordinate of its center of mass and the volume-averaged rising speed. Additionally, the shape of the bubble at $t = 3$ is extracted from the data and compared to a reference solution \cite{hysing2009quantitative}. These results are depicted in Figures \ref{fig:bubble1}, to reasonable agreement with the reference.
    In the second, more challenging setup, we have: $\rho_\mathrm{water} = 1000$, $\rho_\mathrm{air} = 1$, $\mu_\mathrm{water} = 10$, $\mu_\mathrm{air} = 0.1$, $\sigma = 1.96$ and $g = 0.98$. Those results are shown in Figure \ref{fig:bubble2}.

    \section{Shock wave and a water column}
    This qualitative benchmark, taken from \cite{chiocchetti2021high}, invastigates an iteraction of a shock wave in air with a column of water. The domain is a square $\Omega = (-D,D)^2$, where $D = 20\unit{mm}$, which is filled with air, except for a circle of radius $R= 3.2\unit{mm}$ in the center, which is of water. In this test, the initial density of air is $\rho_\mathrm{air} = 1.18 \unit{kg/m^3}$ and the density of water is $\rho_\mathrm{water} = 998.2\unit{kg/m^3}$. The initial pressure of both fluids is $p = 10^5 \unit{Pa}$ (atmospheric pressure). The air is treated as a perfect fluid with $\gamma = 1.4$. A planar shock wave is defined in the left part of the domain and propagates to the right at a speed $v_\mathrm{shock} = 1.3 c_\mathrm{air}$. The simulation starts at $t = -1\unit{\mu s}$ and the initial distance between the shock and the leftmost point of the cylinder is $d = (1\unit{\mu s})\times v_\mathrm{shock}$, such that the impact occurs exactly at $t=0\unit{\mu s}$. The post-shock values, given as a numerical solution to Rankine-Huguniot equations, are prescribed in the initial state to the left from the discontinuity. The surface tension between water and air is $\sigma = 0.072\unit{N/m}$.

\begin{figure}[htb!]
        \centering
        \begin{subfigure}{0.33\linewidth}
            \includegraphics[width=\linewidth]{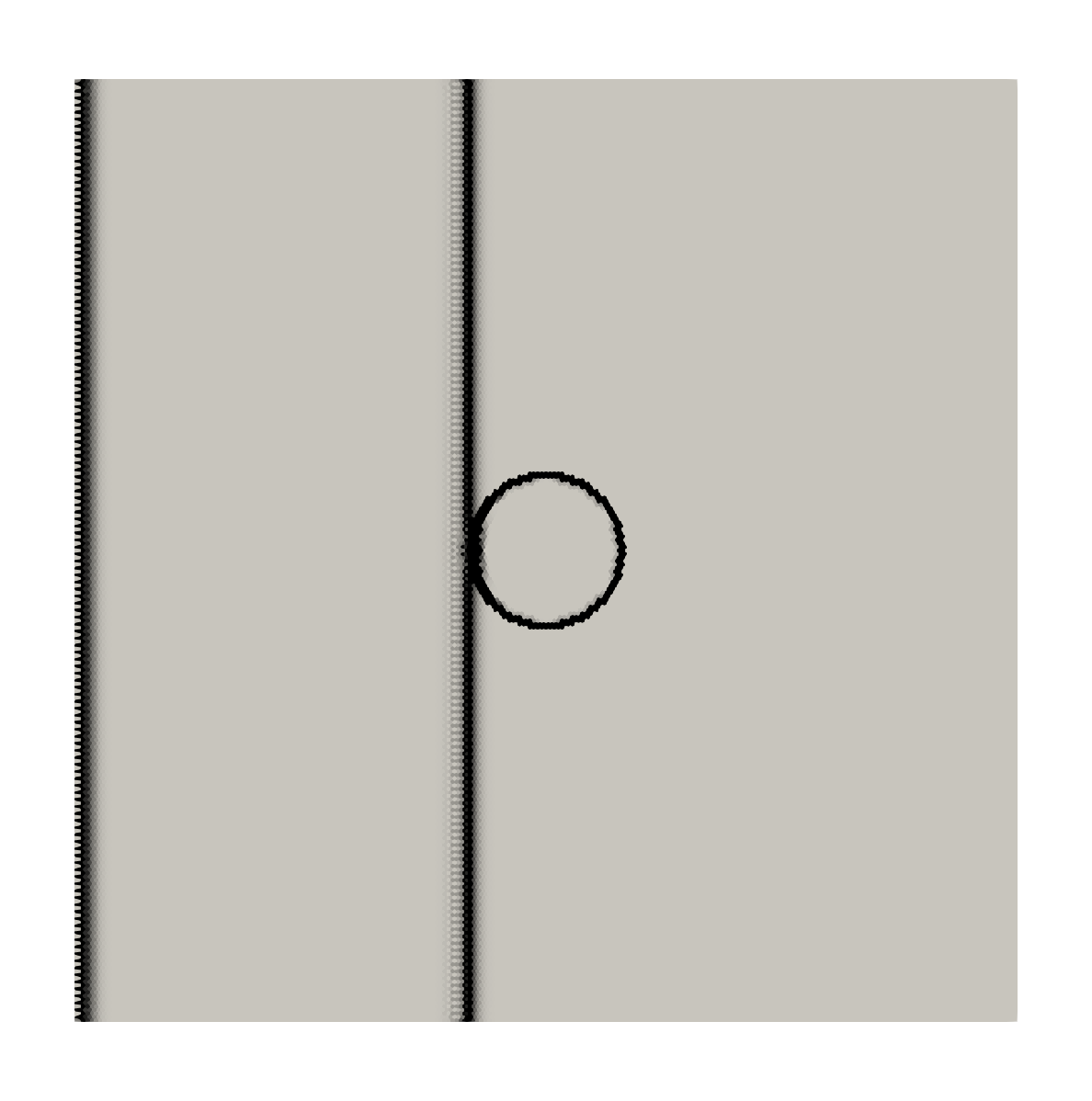} 
            \caption{$t = 0\unit{\mu s}$}
        \end{subfigure}%
        \begin{subfigure}{0.33\linewidth}
            \includegraphics[width=\linewidth]{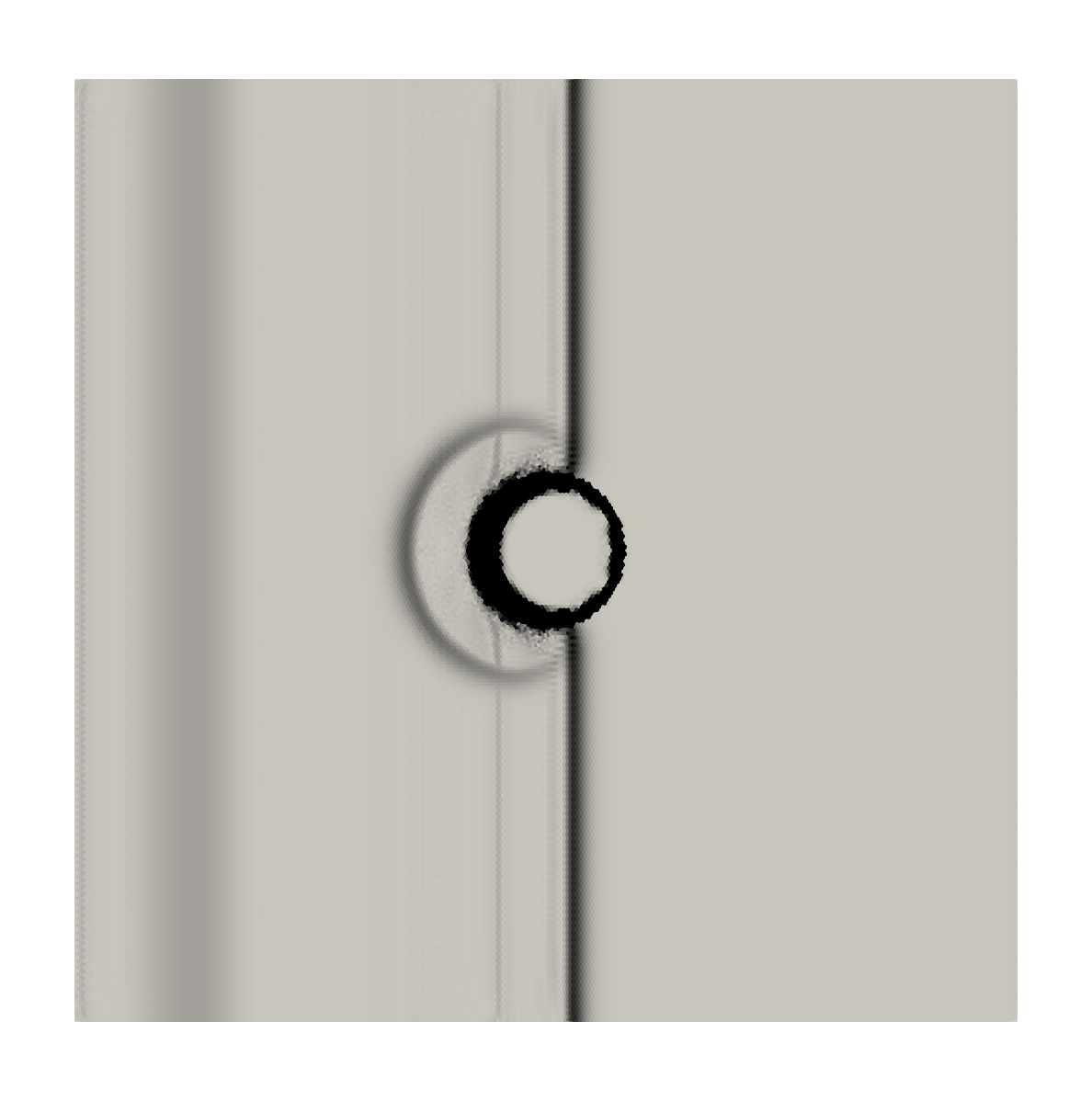} 
             \caption{$t = 10\unit{\mu s}$}
        \end{subfigure}%
        \begin{subfigure}{0.33\linewidth}
            \includegraphics[width=\linewidth]{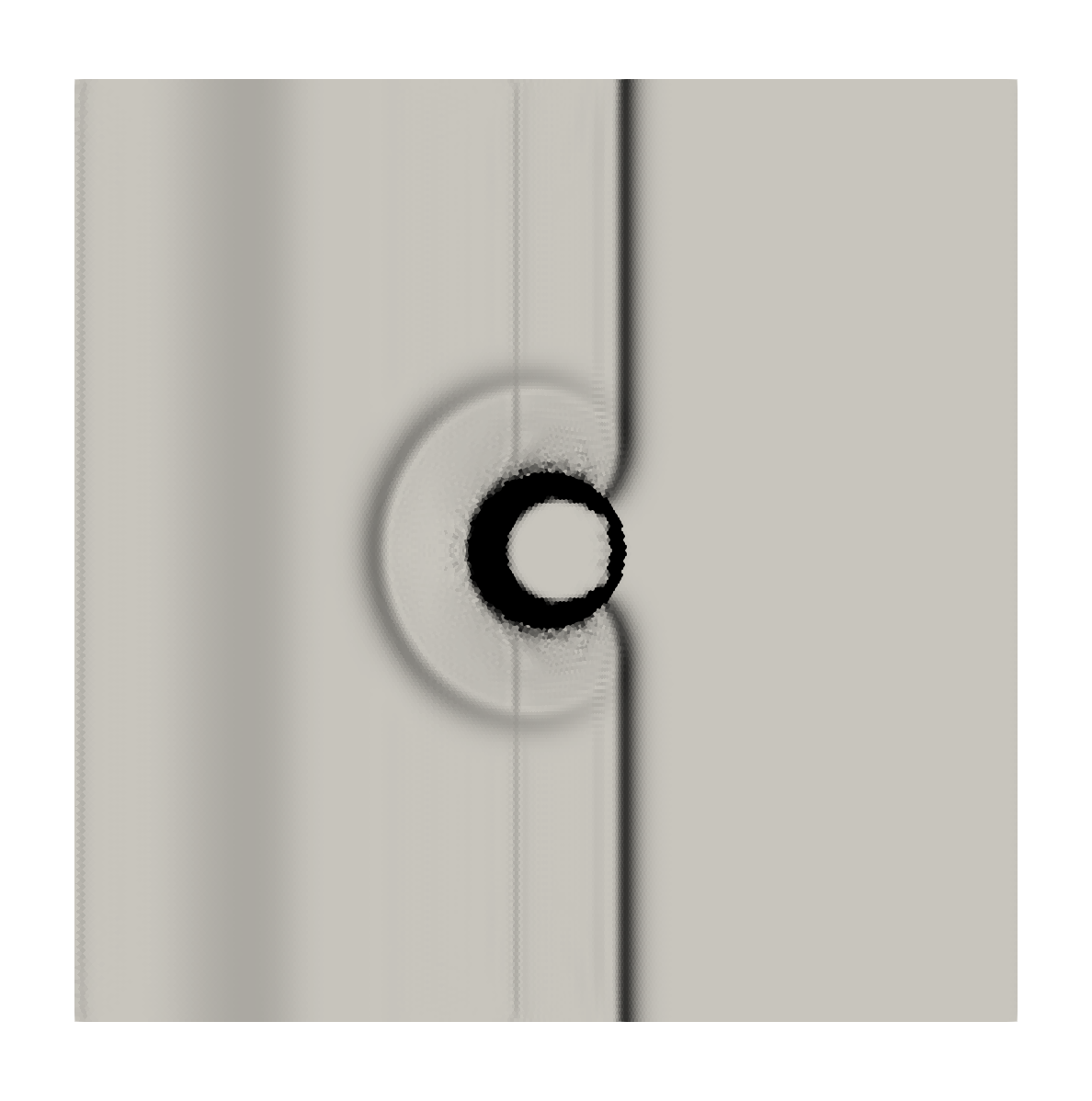} 
             \caption{$t = 15\unit{\mu s}$}
        \end{subfigure}
        \begin{subfigure}{0.33\linewidth}
            \includegraphics[width=\linewidth]{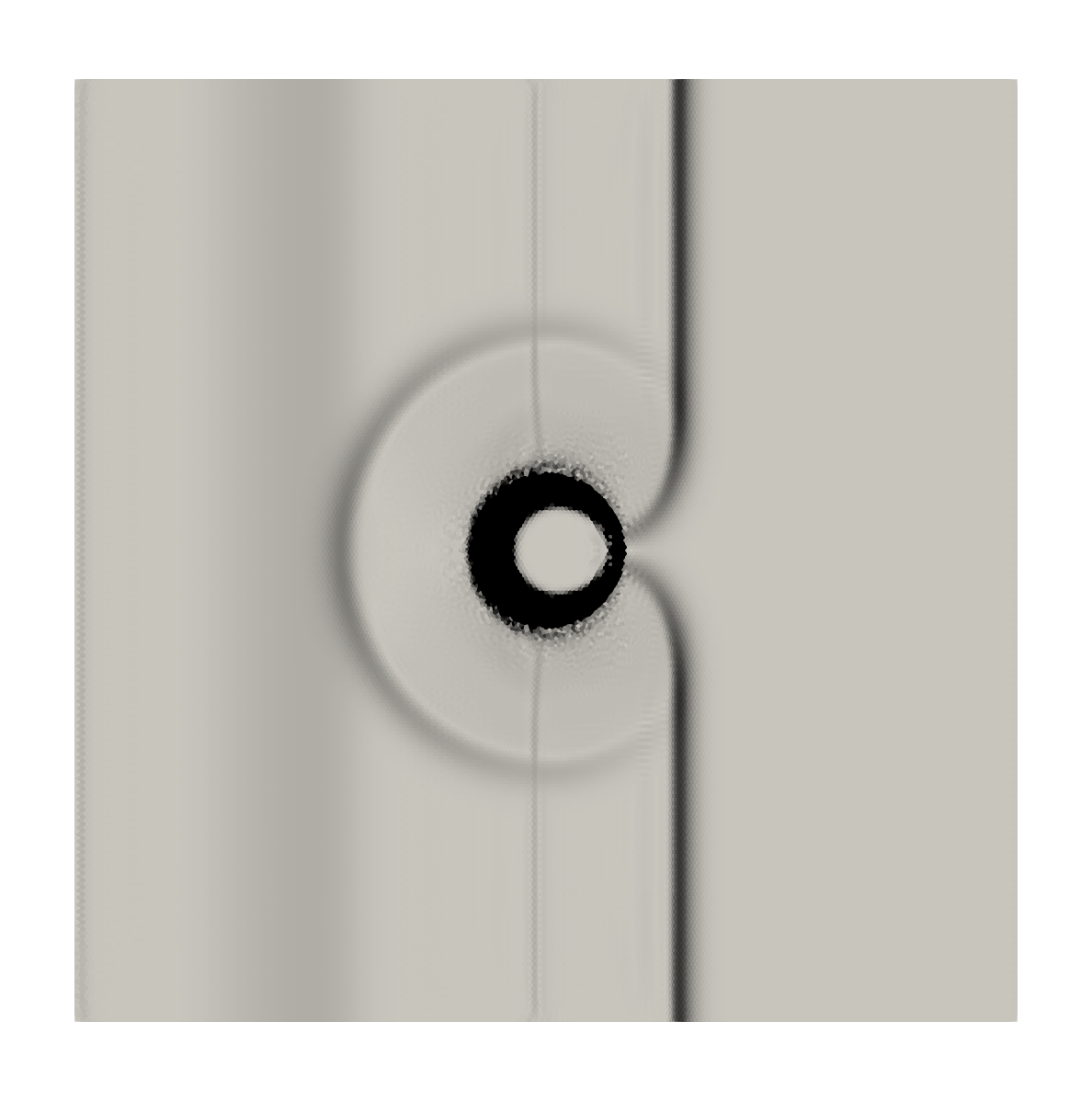} 
             \caption{$t = 20\unit{\mu s}$}
        \end{subfigure}%
        \begin{subfigure}{0.33\linewidth}
            \includegraphics[width=\linewidth]{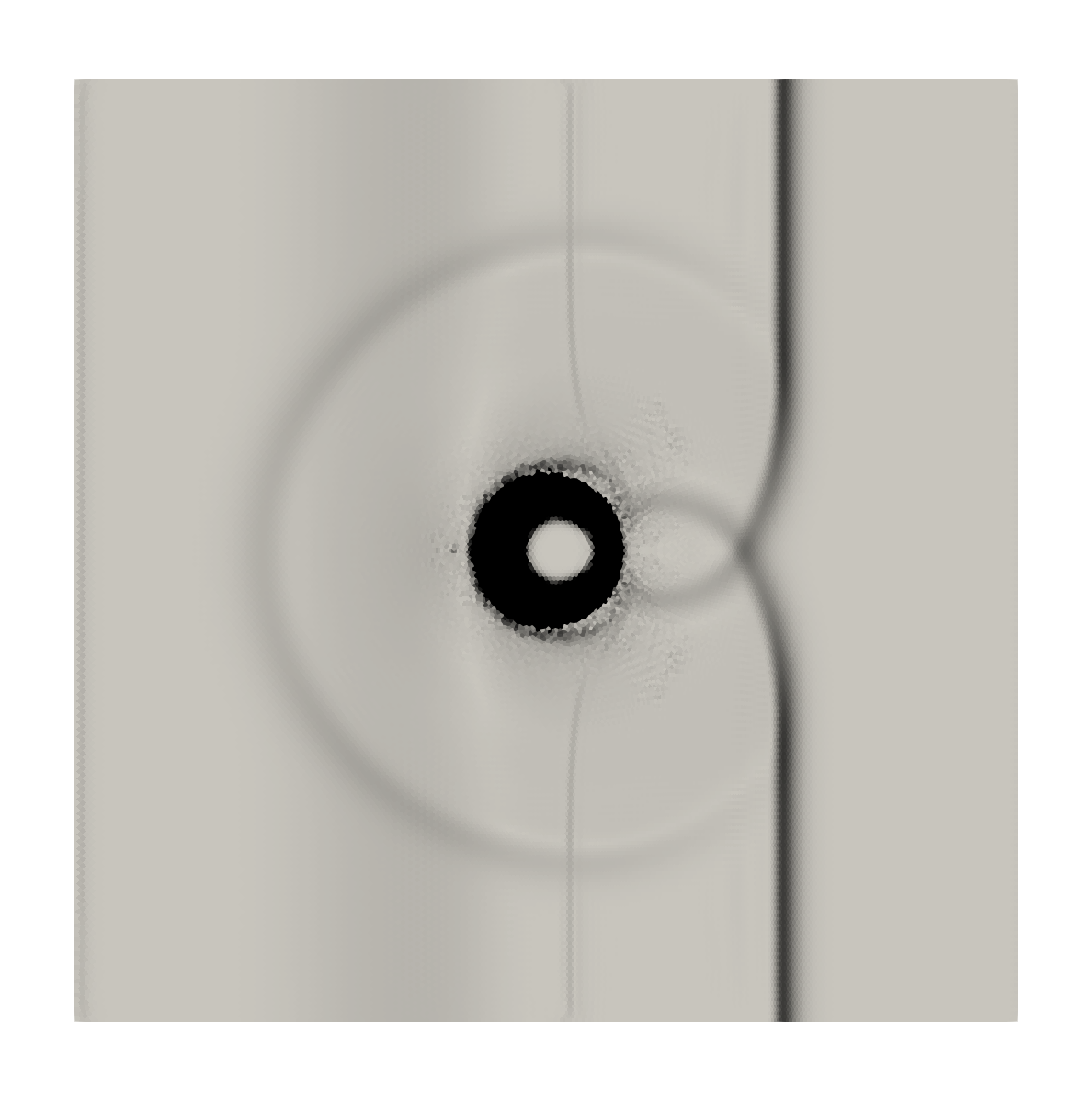}
             \caption{$t = 30\unit{\mu s}$}
        \end{subfigure}%
        \begin{subfigure}{0.33\linewidth}
            \includegraphics[width=\linewidth]{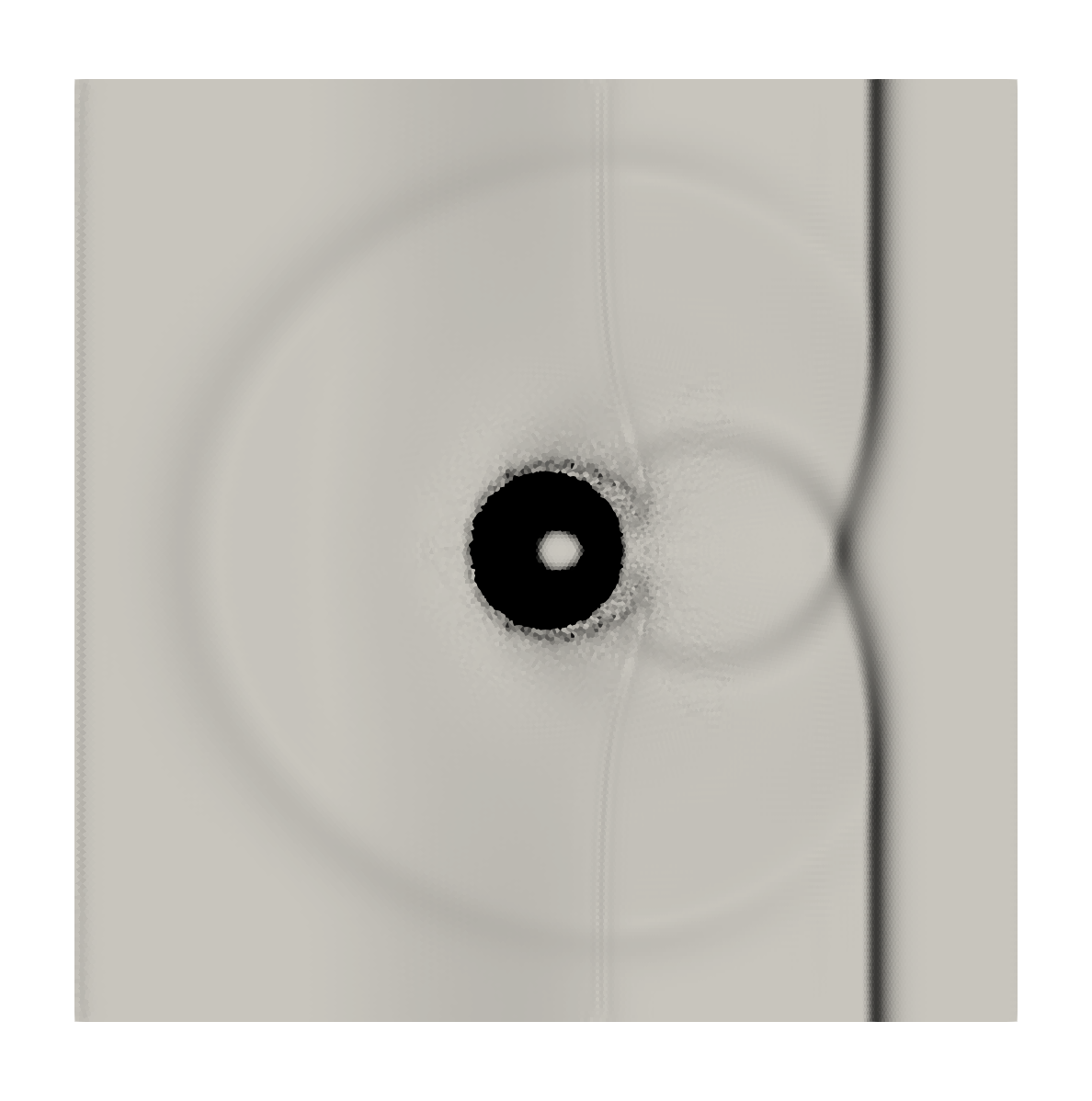} 
             \caption{$t = 40\unit{\mu s}$}
        \end{subfigure}
        \begin{center}
            \includegraphics[width=0.5\linewidth]{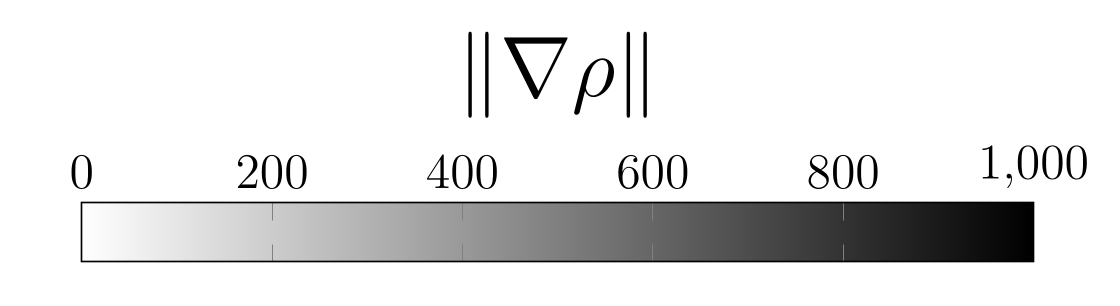}
        \end{center}
        \caption{Numerical schlieren images of the column-shock interaction.}
        \label{fig:shock-column}
    \end{figure}

    The results of the simulation in Figure \ref{fig:shock-column} verifies the ability of our numerical method to capture shock waves and sharp interfaces at the same time.

    \section{Conclusion}
    In this paper, we introduced a technique based on quasi-Lagrangian Voronoi cells for multi-phase flows with sharp interfaces. The presented method is suitable for large density ratios (up to one thousand). The accuracy and the robustness of the method was tested in six benchmark, showing good agreement in comparison to analytical results, reference solutions and experimental data whenever available. The method can handle shock waves as well as surface tension.

    Currently, the two main limitations are the lack of theoretical analysis and that the approximation of derivatives is of first order. Furthermore, our code only works for two-dimensional problems. However, all the underlying ideas can be extended to three dimensional space. We foresee an investigation of the possible extension to higher order of accuracy as well.
    
    \section*{Acknowledgments}
    % This work received financial support by the Italian Ministry of University
    % and Research (MUR) in the framework of the PRIN 2022 project No. 2022N9BM3N
    % (CUP E53D23005840006)
    
    This work was funded by the European Union -- Next Generation EU, Mission 4
    Component 2 -- CUP E53D23005840006. WB acknowledges research funding by
    "Agence Nationale de la Recherche" (ANR) project No. ANR-23-EXMA-0004. WB is
    a member of the INdAM GNCS group in Italy.

\printbibliography
	
\end{document}